\documentclass[11pt]{amsart}

\usepackage{amssymb}
\usepackage{latexsym}
\usepackage{amsmath}

\usepackage{amsthm}

\usepackage{amsfonts}
\usepackage{color}

\usepackage{pdfsync}
\usepackage[usenames,dvipsnames]{pstricks}
\usepackage{epsfig}
\usepackage{pst-grad} 

\usepackage{pst-plot} 

\newcommand{\R}{\mathbb{R}}
\newcommand{\N}{\mathbb{N}}
\newcommand{\Z}{\mathbb{Z}}
\newcommand{\Q}{\mathcal{Q}}

\newcommand{\I}{\mathcal{I}}
\newcommand{\J}{\mathcal{J}}
\newcommand{\T}{\mathbb{T}}
\newcommand{\Tr}{\mathrm{Tr}}
\newcommand{\B}{\mathcal{B}}
\newcommand{\osc}{\mathrm{osc}}

\theoremstyle{plain}
\newtheorem{defi}{Definition}[section]

\newtheorem{prop}[defi]{Proposition}
\newtheorem{teo}[defi]{Theorem}
\newtheorem{cor}[defi]{Corollary}
\newtheorem{lema}[defi]{Lemma}

\theoremstyle{definition}

\newtheorem{remark}{Remark}

\theoremstyle{remark}

\numberwithin{equation}{section}

\definecolor{vert}{rgb}{0,0.4,0}

\begin{document}

\title[]{Lipschitz Regularity for Integro-Differential Equations with Coercive Hamiltonians and Application to Large Time Behavior.}

\author[]{Guy Barles}

\address{
Guy Barles:
Laboratoire de Math\'ematiques et Physique Th\'eorique (UMR CNRS 7350), F\'ed\'eration Denis Poisson (FR CNRS 2964),
Universit\'e Fran\c{c}ois Rabelais Tours, Parc de Grandmont, 37200 Tours, FRANCE. 
\newline {\tt Guy.Barles@lmpt.univ-tours.fr}
}

\author[]{Olivier Ley}
\address{ 
Olivier Ley: IRMAR, INSA de Rennes, 35708 Rennes, FRANCE. 
\newline {\tt  olivier.ley@insa-rennes.fr}
}

\author[]{Erwin Topp}
\address{
Erwin Topp:
Departamento de Matem\'atica y C.C., Universidad de Santiago de Chile,
Casilla 307, Santiago, CHILE.
\newline {\tt erwin.topp@usach.cl}
}

\date{\today}

\begin{abstract} 
In this paper, we provide suitable adaptations of the ``weak version of Bernstein method'' introduced by the first author in 1991, in order to obtain Lipschitz regularity results and Lipschitz estimates for nonlinear integro-differential elliptic and parabolic equations set in the whole space. Our interest is to obtain such Lipschitz results to possibly degenerate equations, or to equations which are indeed ``uniformly elliptic'' (maybe in the nonlocal sense) but which do not satisfy the usual ``growth condition'' on the gradient term allowing to use (for example) the Ishii-Lions' method. We treat the case of a model equation with a superlinear coercivity on the gradient term which has a leading role in the equation. This regularity result together with comparison principle provided for the problem allow to obtain the ergodic large time behavior of the evolution problem in the periodic setting. 
\end{abstract}



\maketitle


\section{Introduction}

The starting point of this article and its main motivation comes from the study of the large time behavior of solutions of nonlinear, nonlocal parabolic partial differential equations. This study requires, in general, two main arguments : Lipschitz estimates which are needed both to prove the compactness of solutions of the evolution equation and to solve the expected stationary limit {\em ergodic problem}; and a Strong Maximum Principle for either the stationary and/or evolution equation to actually prove the convergence. It is worth pointing out that such Strong Maximum Principle is often obtained through a linearization of the equation, which also uses the gradient bound and therefore Lipschitz estimates may be also used indirectly.
In this short description of the method, we would like to stress on the fact that Lipschitz estimates play a central role in all the steps. 

In order to be more specific, we turn to \cite{Barles-Souganidis1} where the large time behavior of solutions of
 (local) nonlinear parabolic PDEs  is studied through two main cases : the {\em sub and super quadratic cases}, the point being that the Lipschitz estimates are obtained in different ways in these two cases. In the subquadratic case, this Lipschitz estimate comes from the uniform elliptic second-order operator (the Laplacian in  \cite{Barles-Souganidis1}) and the subquadratic assumption on the nonlinear terms also related with the $x$-dependence of the Hamiltonian. Technically, this is done using the Ishii-Lions' method \cite{Ishii-Lions}.

On the contrary, in the superquadratic case,  the Lipschitz estimate comes from the nonlinear term through the {\em weak Bernstein's method} (\cite{Barles2}), which has the advantage of being able to handle degenerate cases and Hamiltonians with arbitrary growth (as the classical Bernstein's method).
 
For nonlocal equations, this program is carried out in the ``subquadratic'' case in a series of papers : the Lipschitz estimate is obtained in \cite{Barles-Chasseigne-Imbert-Ciomaga-lip} and the large time behavior in \cite{Barles-Chasseigne-Ciomaga-Imbert} using the Strong Maximum Principle of \cite{Ciomaga}. In the superquadratic case, the contribution of \cite{Barles-Koike-Ley-Topp} is to obtain $C^{0,\alpha}$-type estimates which are sufficient to obtain some large time behavior but for purely nonlocal operator (no mixing of second-order differential operator and nonlocal one). This is one of the rare cases where the Lipschitz estimates can be avoided.

The aim of this paper is to complete this study by providing for some model equations Lipschitz regularity results by using a weak Bernstein's method for bounded viscosity solutions of such nonlocal PDEs.

To the best of our knowledge, there is no general extension of the weak Bernstein's method to the case of nonlocal equations. The reason for such lack of extension may come from the fact that, for PDEs, Bernstein's method (weak or classical) uses a change of variable and such change is not easy to handle for nonlocal equations. But in \cite{Barles-Souganidis1}, only one exponential change is used and it turns out that it fits well with a variety of nonlocal equations including models which are not covered by previous results. To show it, we have decided to treat the case of a rather simple equation, but involving the relevant difficulty, in order to emphasize the new point, namely the additional needed estimates to treat the nonlocal part of the equation.

Our model equations are
\begin{equation}\label{eq}
\lambda u - \Tr(A(x) D^2 u) - \I^j(u,x) + H(x, Du) = 0 \quad \mbox{in} \  \R^d, 
\end{equation}
in the stationary case and its time-dependent version
\begin{equation}\label{pareqIto}
\partial_t u - \mathrm{Tr}(A(x) D^2u) - \I^j(u(\cdot, t), x) + H(x, Du) = 0 \quad \mbox{in} \ Q,
\end{equation}
where $Q = \R^d \times (0,+\infty)$. In both cases, $u : \R^d \to \R$ is the unknown function, $Du, D^2 u$ denote respectively its gradient and Hessian matrix. The main assumptions are $\lambda \geq 0$, $H \in C(\R^d \times \R^d)$ is superlinear in the gradient term, $A$ takes values in the set $\mathbb{S}^{d}_+$ of nonnegative symmetric matrices and $\I^j$ is a nonlocal operator in the \textsl{L\'evy-Ito} form, defined as
\begin{equation}\label{operator}
\I^j(\phi, x) = \int_{\R^d} [\phi(x + j(x,z)) - \phi(x) - \mathbf{1}_B (z)
\langle D\phi(x), j(x, z) \rangle] \nu(dz),
\end{equation}
for $x \in \R^d$. Here $\phi$ is a bounded function which is $C^2$ in a neighborhood of $x$, the function $j : \R^d \times \R^d \to \R^d$ is the \textsl{jump function} and $\nu$ is a \textsl{L\'evy type measure}, which 
is regular and nonnegative. Finally $\mathbf{1}_B$ is the indicator function of the unit ball $B$. Precise assumptions over the data will be given later on.
Note that $-\I^j$ is, up to a normalizing constant, the fractional Laplacian of order $\sigma \in (0,2)$ when
$j(x,z)=z$ and $\nu(dz)=|z|^{-d-\sigma}dz,$ see~\cite{Hitch}.

We emphasize the fact that equations \eqref{eq} and \eqref{pareqIto} may be degenerate in the second-order and/or in the nonlocal term. Hence, our results rely on the coercivity of $H$ in the gradient term, but, in contrast to~\cite{Barles-Koike-Ley-Topp} we can get the result in some cases when the growth degree of coercivity of $H$ in $p$ is less than the order of the nonlocal operator.

We recall that, for local equations, the formal idea of the Bernstein's method is to show that $|Du|^2$ is a subsolution of a suitable elliptic equation. The estimate of $Du$ is then obtained by applying the Maximum Principle, either in a bounded domain or in the whole space.  
The analysis provided in the introduction of~\cite{Barles2} shows that this is possible if the equation satisfies some ``structure condition''. However, such a condition is not directly verifiable at first glance and it is necessary to perform a change of variables which leads to a new equation satisfying this property. The second key information in~\cite{Barles2} is that the formal analysis, consisting in differentiating the equation and therefore requiring smooth solutions, can be justified by viscosity solutions' method and therefore for just continuous solutions.

As in~\cite{Barles2}, in our case the impossibility to differentiate the equation is carried out by the above mentioned viscosity argument that allows us to contrast the Lipschitz bounds of the solution to the problem 
with respect to the Lipschitz bounds of the data, and for this reason we must restrict ourselves to Lipschitz $x$-dependent 
problems.

Of course, the main difference of the application of the method in
the current setting is the presence of the nonlocal term. Recalling the
definition of $\I^j$ in~\eqref{operator}, it is worth to mention that the application of Bernstein
method when $|\nu| < \infty$ and/or $j$ does not depend on $x$ provide
Lipschitz bounds with few extra efforts compared with the already known
second-order case. When the measure is finite there is no differential
effect coming from the nonlocal term in the equation and the
corresponding term is easily controllable, while if $j$ does not depend
on $x$, then the operator is translation invariant which is a favorable
situation in the Bernstein's method where the $x$-dependence of each
differential term is important. For this reason we concentrate in
details on the most difficult scenario of singular measures $\nu$ and
$x$-dependence of $j$ in the definition of $\I^j$, and whose treatment is summarized through Lemma~\ref{estim-nonlocale} below. As in the local
setting, Lipschitz conditions must be requested on $j$, ad-hoc to the
integral configuration of the problem, and these assumptions are
sufficient to control the influence of the nonlocal term with the
stronger coercive effect of the gradient term, no matter the ``order" of
the singularity of the L\'evy measure $\nu$ is.

In the last section of this paper, we apply these Lipschitz regularity results to the study of the large time behavior of the associated evolution problem in the periodic setting. As we already mentioned above, one of the main consequences of this regularity result
is a ``linearization'' procedure of the Hamiltonian that allows us to prove
a version of the Strong Maximum Principle 
for the evolution problem. Roughly speaking, after the mentioned linearization procedure, we can propagate the maximum 
value of a solution of the corresponding linearized problem in the directions of the uniform ellipticity of the 
second-order term as it is performed by Bardi and Da Lio in~\cite{Bardi-DaLio}, 
meanwhile it is propagated in the directions of the degeneracy of the second-order term due to a covering property 
of the support of the measure defining the nonlocal term in the flavour of Coville~\cite{Coville, Coville2}. The novelty is to combine in a better way these two types of (very
different) arguments: this leads to a simpler formulation and a slight
improvement of the results of Ciomaga in~\cite{Ciomaga}.

Once comparison principle, Lipschitz regularity and strong maximum principles are available, 
we follow the lines presented in~\cite{Barles-Souganidis} for first-order 
equations,~\cite{Barles-Souganidis,Tchamba} 
for second-order equations and~\cite{Barles-Chasseigne-Ciomaga-Imbert,Barles-Koike-Ley-Topp} for nonlocal problems
to conclude the solution of the evolution problem behaves, up to a linear factor in time, as 
the solution of the so-called \textsl{ergodic problem}, which can be understood as an homogenization of~\eqref{eq}
when we let $\lambda \to 0$.

We finish this introduction section mentioning that most of the results of this paper can be extended to equations which are  nonlinear in the second-order and nonlocal term, such as Bellman-Isaacs-type nonlinearities arising in game theory. This can be explained by the homogeinity of such operators together with its sub/superaditivity related to Pucci-type associated extremal operators, which do not change the arguments consistently provided there is a weak coupling with the gradient term. However, we do not pursue in this direction for simplicity of the presentation.

\bigskip

\noindent
{\bf Basic notation and organization of the paper.}
In this paper we consider the notion of viscosity solution, see~\cite{Alvarez-Tourin, Barles-Imbert} for 
a definition of this concept in the integro-differential framework.

We use the notation $\mathrm{USC},$ $\mathrm{LSC},$  $\mathrm{C_b}$ and $\mathrm{BUC},$
for upper and lower semicontinuous functions, continuous bounded functions and bounded uniformly continuous
functions, respectively.

We recall
that $Q = \R^d \times (0,+\infty)$, we write $Q_T = \R^d \times (0, T]$ for $T > 0$.
For $a > 0$ and $x \in \R^d$ we denote $B_a(x)$ the ball of center $x$ and radius $a$, $B_a$ when $x$
is the origin and simply $B$ with in addition $a=1$.

For a set $A \subset \R^d$, $x, p \in \R^d$ and $\phi$ a bounded function, we define
\begin{equation*} 
\I^j[A](\phi, x,p) = \int_{A} [\phi(x + j(x,z)) - \phi(x) 
- \mathbf{1}_B (z)\,\langle p, j(x, z) \rangle] \nu(dz).
\end{equation*}

If $A= \R^d$ we write $\I^j(\phi, x,p) := \I^j[\R^d](\phi, x, p).$
When $\phi \in C^2(\R^d)\cap L^\infty(\R^d)$ we also write
$\I^j[A](\phi, x) := \I^j[A](\phi, x, D\phi(x))$.

Note that under mild assumptions on the jump function $j$ and on the measure
$\nu$ (see (M) below), by its definitions the operator~\eqref{operator} is well-defined for a bounded smooth function $\phi$ and can be written as
\begin{equation*}
\I^j(\phi, x)= \I^j[\R^d](\phi, x, D\phi(x)).
\end{equation*}

The paper is organized as follows: in section~\ref{comparisonsection} we provide the well-posedness to the stationary and evolution problem in a general framework. Section~\ref{Lipschitzsection} is devoted to the main result of the paper, which is the Lipschitz regularity for solutions of the stationary and evolution problem
using the weak Bernstein method. Finally, in section~\ref{LTBsection} we restrict ourselves to the periodic setting and provide a strong maximum principle from
which we deduce the large time behavior result.




\section{Comparison principle and its consequences.}
\label{comparisonsection}

In this section we study the well-posedness of the parabolic problem when (\ref{pareqIto}) is associated with the initial data
\begin{equation}
\label{initialdata} u(\cdot,0) = u_0 \quad \mbox{in} \ \R^d,
\end{equation}
where $u_0$ is, at least, a bounded, continuous function in $\R^d$. The initial condition is
satisfied in the sense of viscosity solutions which reduce here to the classical sense by
Lemma~\ref{lemat=0} (see also Remark~\ref{remt0}).

The well-posedness follows by comparison principle for bounded sub and supersolutions obtained
under the following assumptions.

\medskip

\noindent
{\bf (A)} \textsl{There exists a continuous function $\sigma : \R^d \to \R^{d \times k}$, $k \leq d$, such 
that $A(x) = \sigma(x) \sigma^T(x)$ for each $x \in \R^d$ and there exists $L_\sigma \geq 0$ such that
\begin{equation*}
|\sigma(x)|  \leq L_\sigma, \quad 
|\sigma(x) - \sigma(y)|  \leq L_\sigma |x - y|, \quad \mbox{for all} \ x, y \in \R^d.
\end{equation*}
}

\medskip

\noindent
{\bf (H1)} \textsl{There exists $m > 1$ and $K, b_m > 0$ such that for all $\mu \in (0,1)$, $x, p \in \R^d$
\begin{equation*}
\mu H(x, \mu^{-1} p) - H(x, p) \geq (1 - \mu) \Big{(} b_m |p|^m - K \Big{)}.
\end{equation*}
}

\medskip

\noindent
{\bf (H2)} \textsl{Let $m$ be as in (H1). There exist moduli of continuity $\zeta_1, \zeta_2$ such that, for all $x, y, p, q \in \R^d$, $|q|\leq 1$, we have
\begin{eqnarray*}
H(y, p + q) - H(x, p) \leq \zeta_1(|x - y|) (1 + |p|^m) + \zeta_2(|q|)(1+|p|^{m - 1}).
\end{eqnarray*}
}

\medskip

\noindent
{\bf (M)} \textsl{ There exists $C_{\nu,j} > 0$ such that
\begin{equation*}
\int \limits_{B^c} \nu(dz) \leq C_{\nu,j}
\quad and \quad
\int \limits_{B} |j(x,z)|^2 \nu(dz) \leq C_{\nu,j},
\end{equation*}
where $a\wedge b=\min(a,b)$.
}
\medskip

\noindent
{\bf (MJ)} \textsl{For any $R \geq 1$, there exists constant $C^0(R), C^1(R), C^2(R) > 0$ such that for all $x, y \in \R^d$
$$
\int_{B_R\setminus B} |j(x, z)| \nu(dz) \leq C^0(R) \; ,
$$
$$
\int_{B_R\setminus B} |j(x, z) - j(y, z)| \nu(dz) \leq C^1(R)|x - y| \; ,
$$
$$
 \int_{B_R} |j(x, z) - j(y, z)|^2 \nu(dz) \leq C^2(R)|x - y|^2\;. 
$$}

We point out that (A) is a classical assumption to prove uniqueness for degenerate equations, see~\cite{Ishii-Lions}.

The hypothesis on the nonlocal term are classical for L\'evy-It\^o operators. Assumption (M) is the so-called \textsl{L\'evy condition} over $\nu$, and it allows to give a sense to the nonlocal operator for bounded $C^2$ functions. On the other hand, (MJ) is a continuity condition to treat the nonlocal terms in the comparison proof. 

Concerning the conditions over the Hamiltonian, 
(H1) gives the structure allowing to apply the weak Bernstein method.
When $H$ is smooth, (H1) reduces to $H_p(x,p)p-H(x,p)\geq b_m|p|^m-K$ (see~\cite{Barles-Souganidis, Ley-Vihn}).
In particular, $H$ is coercive, see~\eqref{Hcoercive}. Assumption
(H2) states the continuity of the Hamiltonian, 
which is relative to its degree of coercivity $m$ stated in (H1). Examples of such Hamiltonians can be found 
in~\cite{Barles-Koike-Ley-Topp} and references therein (see also~\cite{GKLV} for some examples described in a detailed way).

Notice that (H2) implies the function $x \mapsto H(x,0)$ is uniformly continuous.


\subsection{Comparison principle.}


\begin{prop}\label{parcomparison}
Let $u_0 \in C_b(\R^d)$, $A$ satisfying (A), $\I^j$ defined as in~\eqref{operator} satisfying (M) and (MJ). Assume $H$ satisfies (H1),(H2). 
Let $u$ be an USC subsolution and $v$ a LSC supersolution to
the problem~\eqref{pareqIto}-\eqref{initialdata} such that $u, v$ are bounded in $\R^d \times [0,T]$ for each $T > 0$. 
Then, 
\begin{equation*}
u \leq v \quad \mbox{in} \ \bar{Q}. 
\end{equation*}
\end{prop}


Before proving the above proposition we introduce some technical 
results that are going to be used in different frameworks in the rest of the paper.
The cell tool is the power function
\begin{equation}\label{varphi}
\varphi_\alpha(x, y) = |x - y|^\alpha, \quad \mbox{for} \ x, y \in \R^d,
\end{equation}
where $ \alpha > 0$. This function will play different roles on the arguments to come, taking into account different values of $\alpha$. 

For $x, y \in \R^d$ with $x \neq y$, define the matrix 
\begin{equation}\label{Zalpha}
Z_\alpha = I_d + (\alpha - 2)(\widehat{x - y}) \otimes (\widehat{x - y}), 
\end{equation}
where $I_d$ is the $d \times d$ identity matrix and $\hat{x} = x/|x|$ for $x \neq 0$. 

With this notation, note that for $x \neq y$, direct computations show that
\begin{equation}\label{calcul-deriv134}
\begin{split}
& D_x \varphi_\alpha(x,y) = -D_y \varphi_\alpha(x,y) = \alpha |x - y|^{\alpha - 1} (\widehat{x - y}), \\
& D^2_{xx} \varphi_\alpha(x,y) = D^2_{yy} \varphi_\alpha(x,y) = \alpha |x - y|^{\alpha - 2} Z_\alpha,\\
& D^2_{xy} \varphi_\alpha(x,y) = D^2_{yx} \varphi_\alpha(x,y) = -\alpha |x - y|^{\alpha - 2} Z_\alpha.
\end{split}
\end{equation}

\begin{lema}\label{lemaADvarphi}

Let $\alpha > 0$, $\varphi_\alpha$ defined in~\eqref{varphi} and let $A: \R^d \to \mathbb{S}^d$ satisfying assumption (A).
Let $\bar{x}, \bar{y} \in \R^d$ (with $\bar{x} \neq \bar{y}$ if $\alpha < 2$) and assume there exist two 
matrices $X, Y \in \mathbb{S}^d$ satisfying the inequality
\begin{equation}\label{matrixineq}
\left [ \begin{array}{cc} X & 0 \\ 0 & -Y \end{array} \right ] \leq D_{(x,y)}^2\varphi_\alpha(\bar{x}, \bar{y}).
\end{equation}

Then, we have the estimate
\begin{equation*}
\mathrm{Tr}(A(\bar{x})X - A(\bar{y})Y) \leq \alpha (1 + |\alpha - 2|)L_\sigma^2 |\bar{x} - \bar{y}|^\alpha.
\end{equation*}
\end{lema}

\noindent
{\bf \textit{Proof:}} Multiplying the left inequality~\eqref{matrixineq} by the nonnegative matrix
\begin{equation}\label{matrixIJ}
\left [ \begin{array}{cc}  \sigma(\bar{x}) \sigma^T(\bar{x}) &  \sigma(\bar{x}) \sigma^T(\bar{y}) \\
 \sigma(\bar{y}) \sigma^T(\bar{x}) & \sigma(\bar{y}) \sigma^T(\bar{y}) \end{array} \right ] 
\end{equation}
and taking traces, the resulting inequality drives us to
\begin{equation*}
\begin{split}
& \Tr(A(\bar{x}) X - A(\bar{y}) Y) \\
\leq & \ \Tr \Big{(} \sigma(\bar{x}) \sigma^T(\bar{x}) D^2_{xx}\varphi_\alpha(\bar{x}, \bar{y}) 
+ \sigma(\bar{x}) \sigma^T(\bar{y}) D_{xy}^2 \varphi_\alpha(\bar{x}, \bar{y}) \\ 
& \qquad +  \sigma(\bar{y}) \sigma^T(\bar{x}) D_{yx}^2 \varphi_\alpha(\bar{x}, \bar{y})
+ \sigma(\bar{y}) \sigma^T(\bar{y}) D_{yy}^2 \varphi_\alpha(\bar{x}, \bar{y}) \Big{)}.
\end{split}
\end{equation*}

Using the computations for the derivatives of $\varphi_\alpha$ and the definition of $Z_\alpha$, we have
\begin{equation*}
\begin{split}
& \Tr(A(\bar{x}) X - A(\bar{y}) Y) \\ 
\leq & \ \alpha |\bar{x} - \bar{y}|^{\alpha -2} \Tr\Big{(}(\sigma(\bar{x}) \sigma^T(\bar{x}) - \sigma(\bar{x}) \sigma^T(\bar{y})
- \sigma(\bar{y}) \sigma^T(\bar{x}) + \sigma(\bar{y}) \sigma^T(\bar{y}))Z_\alpha \Big{)} \\
= & \ \alpha|\bar{x} - \bar{y}|^{\alpha - 2} \Tr\Big{(} (\sigma(\bar{x}) - \sigma(\bar{y}))(\sigma(\bar{x}) - \sigma(\bar{y}))^T Z_\alpha\Big{)},
\end{split}
\end{equation*}
but using Schwarz inequality, we obtain that
\begin{equation*}
\Tr(A(\bar{x}) X - A(\bar{y}) Y) 
\leq \alpha (1 + |\alpha- 2|)|\bar{x} - \bar{y}|^{\alpha - 2} |\sigma(\bar{x}) -\sigma(\bar{y})|^2. 
\end{equation*}

Finally, applying condition (A) we conclude the result.
\qed

\medskip

Due to the lack of compactness of $\R^d$, some localization argument in $x$ is needed. In the sequel we use a nonnegative function $\psi \in C^2_b(\R^d)$  satisfying the following properties
\begin{equation}\label{defpsi}
\left \{ \begin{array}{l} \psi = 0 \quad \mbox{in} \ B, \\
\psi = \Psi \quad \mbox{in} \ B^c_2 \quad \mbox{for some constant} \ \Psi >0, \\
0 \leq \psi \leq \Psi \quad \mbox{in} \  B_2 \setminus B; \quad \mbox{and} \\ 
||D\psi||_\infty, ||D^2\psi||_\infty \leq \Lambda \quad \mbox{for some} \ \Lambda > 0. 
\end{array} \right .
\end{equation}

Next lemma states the estimates for a localization  function based on $\psi$.


\begin{lema}\label{lemapsibeta}
Assume (M) hold. Let $\psi$ satisfying the properties listed in~\eqref{defpsi} and for $\beta > 0$, define 
the function 
\begin{equation}\label{psibeta}
\psi_\beta(x) = \psi(\beta x), \quad x \in \R^d. 
\end{equation}

Then, $\psi_\beta$ satisfies
\begin{eqnarray*}
&& ||D\psi_\beta||_\infty \leq  \Lambda \beta, \quad ||D^2 \psi_\beta||_\infty \leq \Lambda \beta^2, \\
&&  ||\I^j[B_\delta\cap A](\psi_\beta, \cdot)||_\infty \leq \Lambda \beta^2 o_\delta(1),\\
&& ||\I^j[B_\delta^c \cap A](\psi_\beta, \cdot)||_\infty
\leq \Lambda o_\beta(1), 
\end{eqnarray*}
where $o_\beta(1), o_\delta(1) \to 0$ as $\beta, \delta \to 0$
respectively and $o_\beta(1)$ depends only on $\Psi,$ $\nu,$ $j$ and
$o_\delta(1)$ depends only on $\nu,$ $j.$
\end{lema}


\noindent
{\bf \textit{Proof of Lemma~\ref{lemapsibeta}.}}
The estimates for $D\psi_\beta$ and $D^2 \psi_\beta$ are obvious. 

Now we consider $\beta > 0$, $0<\delta \leq 1,$ $x\in\R^d$ and $A\subset \R^d$ measurable. Using the smoothness of $\psi,$ we have
\begin{eqnarray*}
\I^j[B_\delta\cap A](\psi_\beta, x)
= \frac{1}{2} \int_{B_\delta \cap A}\int_0^1 
\langle D^2\psi_\beta(x+\theta j(x,z))j(x,z),j(x,z)\rangle d \theta \nu(dz),
\end{eqnarray*}
from which, we easily deduce that
\begin{equation*}\label{estim1111}
|\I^j[B_\delta\cap A](\psi_\beta, x)|\leq 
\frac{1}{2}||D^2\psi_\beta||_{\infty}
\int_{B_\delta\cap A} |j(x, z)|^2 \nu(dz). 
\end{equation*}

Then, using (M) and the estimates for $D^2 \psi_\beta$ we get
\begin{eqnarray}\label{estim472}
|\I^j[B_\delta\cap A](\psi_\beta , x)|
\leq 
\frac{1}{2}\beta^2\Lambda
\int_{B_\delta\cap A} |j(x,z)|^2 \nu(dz)= \beta^2\Lambda o_\delta (1). 
\end{eqnarray}

From the definition of $\psi_\beta,$ we have, for all $x,y\in\R^d,$
$|\psi_\beta(x)-\psi_\beta(y)|\leq \Psi$ and $\psi_\beta(x)\to 0$ as $\beta\to 0.$
Therefore, since $\int_{B^c}\nu(dz)<\infty$ by (M), from the Dominated Convergence
Theorem, we obtain
\begin{eqnarray}\label{estim473}
\I^j[B^c\cap A](\psi_\beta , x)=o_\beta (1),
\end{eqnarray}
where $o_\beta (1)$ depends only on $\Psi,$ $\nu,$ $j.$

From~\eqref{estim472} (with $\delta=1$) and~\eqref{estim473} it follows
\begin{eqnarray*}
\I^j[A](\psi_\beta , x)=o_\beta (1)+\beta^2\Lambda o_1(1)= o_\beta (1),
\end{eqnarray*}
and from here we finally get that
\begin{eqnarray*}
\I^j[B_\delta^c \cap A](\psi_\beta , \cdot)=\I^j[A](\psi_\beta , \cdot)
- \I^j[B_\delta\cap A](\psi_\beta , \cdot)
= o_\beta (1),
\end{eqnarray*}
from which the result follows.
\qed

\medskip


Next, we have the following
\begin{lema}\label{lemat=0}
Let $A$ satisfying (A), $\I^j$ defined in~\eqref{operator} such that its components satisfy (M), (MJ) and let $H$ satisfying (H1). Let 
$u$ be an USC subsolution and $v$ a LSC supersolution
to problem~\eqref{pareqIto}-\eqref{initialdata}, 
bounded in $\R^d \times [0,T]$ for each $T > 0$. Then, 
$u(x,0) \leq u_0(x) \leq v(x,0)$ for all $x \in \R^d$.
\end{lema}

\begin{remark}\label{remt0}
The initial condition has to be understood in the viscosity sense, see
\cite[Definition 2.1]{Barles-Topp}. Lemma~\ref{lemat=0} states that actually it holds in the classical
sense. We refer to~\cite{Barles-Topp} (see also~\cite{DaLio}) for a proof of this result in the case of Dirichlet problem in bounded domains. In the current setting, the proof must be slightly modified by a standard localization procedure that allows to deal with the lack of compactness of $\R^d$. 
Hence, we remark that uniform continuity of the initial data is not necessary in the
proof of the comparison principle.
\end{remark}

\noindent
{\bf \textit{Proof of Proposition~\ref{parcomparison}}.}
We will argue over the finite horizon problem
\begin{equation*}
\left \{ \begin{array}{rll} \partial_t u - \mathrm{Tr}(A D^2 u) 
- \I^j(u, x) + H(x, Du) & = 0 \quad & \mbox{in} \ Q_T \\ 
u(x, 0) & = u_0(x) \quad & x \in \R^d, \end{array} \right .  
\end{equation*}

We are going to prove that $M:=\sup_{Q_T}\,(u - v) \leq 0$; the general result on $Q$ follows since $T$ is arbitrary.

We argue by contradiction, assuming that $M > 0$ and for $\mu, \eta \in (0,1)$ to be fixed later, we define $(x,t) \mapsto \bar{u}(x,t) := \mu u(x, t) - \eta t$. Since $u$ is bounded, we have $\sup_{Q_T}\,(\bar{u} - v) \geq M/2 >0$ if $\eta$ is small enough and $\mu$ sufficiently close to $1$.
We notice that since $u$ is a subsolution to~\eqref{pareqIto} then $\bar{u}$ satisfies
\begin{equation*} 
\partial_t \bar{u} - \mathrm{Tr}(A D^2 \bar{u}) - \I^j(\bar{u}, x) + \mu H(x, \mu^{-1} D\bar{u}) \leq - \eta \quad \mbox{in} \ Q_T, 
\end{equation*}
in the viscosity sense.

Next we consider for $\epsilon > 0$ the function
\begin{equation*}
(x,y,s,t) \mapsto \Phi(x,y,t) := \bar{u}(x,t) - v(y,t) - \phi(x,y),
\end{equation*}
where $\phi(x,y) := \epsilon^{-2}|x - y|^2 + \psi_\beta(y)$ and $\psi_\beta$ is defined as in Lemma~\ref{lemapsibeta} with a localization function $\psi$ defined as in~\eqref{defpsi} with $\Psi = 2 (||u||_{L^\infty(Q_T)} + ||v||_{L^\infty(Q_T)})$.
With this we see that for all $\beta$ small enough
\begin{equation*} 
\overline{M} := \sup \limits_{\bar{Q}_T \times \bar{Q}_T} \Phi  > 0,
\end{equation*}
and this supremum is achieved at some point $(\bar{x}, \bar{y}, \bar{t}) \in \R^N \times \R^N \times [0,T]$. However, in view of Lemma~\ref{lemat=0} and the positiveness of $\overline{M}$ we necessarily have that $\bar{t} >0$, and standard arguments in the viscosity theory imply this sequence of points (depending on $\epsilon$ and $\beta$) satisfies
\begin{equation}\label{comparison1}
\begin{split}
& \epsilon^{-2}|\bar{x} - \bar{y}|^2 \to 0 \quad\hbox{if  }\epsilon \to 0,\ \beta \ \hbox{fixed,} \\
& \overline{M} \to \sup_{Q_T}\,(\bar{u} - v-\psi_\beta ) \quad\hbox{if  }\epsilon \to 0,\ \beta \ \hbox{fixed,}\\
&\sup_{Q_T}\,(\bar{u} - v-\psi_\beta ) \to \sup_{Q_T}\,(\bar{u} - v) >0 \quad\hbox{when  }\beta \to 0.
\end{split}
\end{equation}

Then we can apply the nonlocal parabolic version of Ishii-Jensen Lemma provided in~\cite{Barles-Imbert} (see also~\cite{Jensen, Ishii3, usersguide} for second-order parabolic equations) which avoids the doubling in the time variable in the definition of $\Phi$ above. Thus, for all $\delta > 0$ and all $\rho > 0$ small enough we have
\begin{equation}\label{testinguvcomparison}
\begin{split}
\varpi - \Tr(A(\bar{x}) X_\rho)  - \I^j[B_\delta](\phi(\cdot, \bar{s}, \bar{y}, \bar{t}), \bar{x}) & \\
 - \I^j[B_\delta^c](\bar{u}, \bar{x}, \bar{p})
+ \mu H(\bar{x}, \mu^{-1}\bar{p}) & \leq o_\rho(1) - \eta, \\
\varpi - \Tr(A(\bar{y}) Y_\rho)  - \I^j[B_\delta](-\phi(\bar{x}, \bar{s}, \cdot, \bar{t}), \bar{y}) & \\
 - \I^j[B_\delta^c](v, \bar{y}, \bar{p} + \bar{q})
+ H(\bar{y}, \bar{p} + \bar{q}) & \geq o_\rho(1)
\end{split}
\end{equation}
where $\varpi \in \R$, $\bar{p} = 2\epsilon^{-2}(\bar{x} - \bar{y}),$ $\bar{q} = -D\psi_\beta(\bar{y})$ and the matrices $X_\rho, Y_\rho \in \mathbb{S}^{d}$ satisfy the inequality
\begin{equation}\label{IJineqcomparison}
-\rho^{-1} I_{2d} \leq \left [ \begin{array}{cc} X_\rho & 0 \\ 0 & -Y_\rho \end{array} \right ] \leq D_{(x,y)}^2\phi(\bar{x}, \bar{s}, \bar{y}, \bar{t}) 
+ o_\rho(1).
\end{equation}

We remark that all the terms $o_\rho(1)$ arising in~\eqref{testinguvcomparison} and~\eqref{IJineqcomparison} 
satisfy $o_\rho(1) \to 0$ for $\epsilon, \beta > 0$ fixed. Subtracting both inequalities in~\eqref{testinguvcomparison} we get
\begin{equation}\label{testingcomparison}
\mathcal{H} \leq \mathcal{A} + \B_\delta + \B^\delta - \eta + o_\rho(1), 
\end{equation}
where
\begin{equation*}
\begin{split}
\mathcal{A} & = \Tr(A(\bar{x}) X_\rho) - \Tr(A(\bar{y}) Y_\rho) \\
\mathcal{H} & = \mu H(\bar{x}, \mu^{-1} \bar{p}) - H(\bar{y}, \bar{p} + \bar{q}) \\
\B_\delta & = \I^j[B_\delta](\phi(\cdot, \bar{s}, \bar{y}, \bar{t}), \bar{x}) 
- \I^j[B_\delta](\phi(\bar{x}, \bar{s}, \cdot, \bar{t}), \bar{y}), \\
\B^\delta & = \I^j[B_\delta^c](\bar{u}(\cdot, \bar{s}), \bar{x}, \bar{p}) - \I^j[B_\delta^c](v(\cdot, \bar{t}), \bar{y}, \bar{p} + \bar{q}).
\end{split}
\end{equation*}

In what follows we estimate each term arising in~\eqref{testingcomparison}. 

\medskip

\noindent
1.- \textsl{Estimate of $\mathcal{A}$:} In view of the definition of $\varphi$ in~\eqref{varphi} we can write
$$
\phi(x,y,t) = \epsilon^{-2} \varphi_2(x,y) + \psi_\beta(y).
$$

Hence, using the estimate given by Lemma~\ref{lemaADvarphi} and applying the estimates for the second derivatives of $\psi_\beta$ in Lemma~\ref{lemapsibeta}, we conclude that
\begin{equation*}
\begin{split}
\Tr(A(\bar{x}) X_\rho - A(\bar{y}) Y_\rho) \leq  2L_\sigma^2 \epsilon^{-2}|\bar{x} - \bar{y}|^{2} + \beta^2 L_\sigma^2 \Lambda  + o_\rho(1).
\end{split}
\end{equation*}

From this, using~\eqref{comparison1}, we conclude
\begin{equation}\label{Acomparison}
\mathcal{A} \leq m_\beta(\epsilon) + o_\beta(1) + o_\rho(1)\,,
\end{equation}
where, here and below, $m_\beta(\epsilon)$ denotes various quantities which tend to $0$ when $\epsilon \to 0$, $\beta$ remaining fixed and $o_\beta(1)\to 0$ as $\beta\to 0,$ $\mu$ and $\eta$ being fixed.

\medskip

\noindent
2.- \textsl{Estimate for $\mathcal{H}$:} Recalling that $|\bar{q}| \to 0$ as $\beta \to 0,$ we have
$|\bar{q}|\leq 1$ for $\beta$ small enough and we can apply (H1) and (H2) to get
\begin{eqnarray*}
\mathcal{H} &\geq& \left[(1 - \mu) b_m - \zeta_1(|\bar{x} - \bar{y}|)\right]|\bar{p}|^m  - \zeta_2(|\bar{q}|)
(1+|\bar{p}|^{m-1}) \\
&& - K(1 - \mu) - \zeta_1(|\bar{x} - \bar{y}|),
\end{eqnarray*}
and applying Young's inequality, we can write
\begin{equation*}
\begin{split}
\mathcal{H} \geq & \ \Big{[}(1 - \mu) b_m - \zeta_1(|\bar{x} - \bar{y}|) - (m - 1)m^{-1} \zeta_2(|\bar{q}|)^{m(m - 1)^{-1}/2} \Big{]}|\bar{p}|^m  \\
& \ - m^{-1}\zeta_2(|\bar{q}|)^{m/2} - \zeta_2(|\bar{q}|)
- K(1 - \mu) - \zeta_1(|\bar{x} - \bar{y}|).  
\end{split}
\end{equation*}

At this point, we fix $\mu =\mu_\eta < 1$ close to 1 in order to have 
$K(1 - \mu_\eta) \leq \eta/4$.
Considering $\epsilon$ and $\beta$ small enough depending on $\mu$ and $b_m$, we can suppose that
$|\bar{x} - \bar{y}|$ and $|\bar{q}|$ are small enough to make positive the term in 
the squared brackets in the last inequality. From this, we conclude that
\begin{equation}\label{Hcomparison}
\mathcal{H} \geq -o_\beta(1) - m_\beta(\epsilon) - \eta/4.
\end{equation}

\medskip

\noindent
3.- \textsl{Estimate for $\B_\delta$:}
Using the definition of $\phi$ and Lemma~\ref{lemapsibeta}, we obtain
\begin{eqnarray*}
\B_\delta &\leq& 
\frac{1}{2}
\frac{||D^2\varphi_2 (\cdot,\bar{y})||_{\infty}}{\epsilon^2}
\int_{B_\delta} |j(\bar{x}, z)|^2 \nu(dz) 
\\
&& \hspace{1.5cm}
+
\frac{1}{2}\left(
\frac{||D^2\varphi_2 (\bar{x},\cdot)||_{\infty}}{\epsilon^2}
+ ||D^2\psi_\beta||_{\infty}
\right)
\int_{B_\delta} |j(\bar{y}, z)|^2 \nu(dz),
\end{eqnarray*}
and therefore, using (M) we conclude that
\begin{equation}\label{B_deltacomparison}
\mathcal{B}_\delta \leq  ( \epsilon^{-2} + \Lambda \beta^2 )o_\delta(1),
\end{equation}
where $o_\delta(1) \to 0$ as $\delta \to 0$.

\medskip

\noindent
4.- \textsl{Estimate for $\B^\delta$:}
We start writing $\B^\delta = \B^{\delta, 1} + \B^{\delta, 2}$ with
\begin{equation*}
\begin{split}
\B^{\delta, 1} & = \I^j[B \cap B_\delta^c](\bar{u}(\cdot, \bar{t}), \bar{x}, \bar{p}) 
- \I^j[B \cap B_\delta^c](v(\cdot, \bar t), \bar{y}, \bar{p} + \bar{q}), \\
\B^{\delta, 2} & = \I^j[B^c](\bar{u}(\cdot, \bar{t}), \bar{x}, \bar{p}) 
- \I^j[B^c](v(\cdot, \bar t), \bar{y}, \bar{p} + \bar{q}).
\end{split}
\end{equation*}

Since $(\bar{x}, \bar{y}, \bar{t})$ is a maximum for $\Phi$, for each $\xi, \xi' \in \R^d$
we have
\begin{equation*}
\begin{split}
& \bar{u}(\bar{x} + \xi, \bar{t}) - v(\bar{y} + \xi', \bar{t}) - \epsilon^{-2}|\bar{x} + \xi - \bar{y} - \xi'|^2 - \psi_\beta(\bar{y} + \xi') \\
\leq & \ \bar{u}(\bar{x}, \bar{t}) - v(\bar{y}, \bar{t}) - \epsilon^{-2}|\bar{x} - \bar{y}|^2 - \psi_\beta(\bar{y}).
\end{split}
\end{equation*}

Using this inequality with $\xi = j(\bar{x}, z), \xi' = j(\bar{y}, z)$ we obtain
\begin{equation}\label{ineg-max127}
\begin{split}
& \bar{u}(\bar{x} + j(\bar{x}, z), \bar{t}) - \bar{u}(\bar{x}, \bar{t}) - (v(\bar{y} + j(\bar{y}, z), \bar{t}) - v(\bar{y}, \bar{t})) \\
\leq & \ \epsilon^{-2}|j(\bar{x}, z) - j(\bar{y},z)|^2 + 2\epsilon^{-2} \langle \bar{x} - \bar{y}, j(\bar{x}, z) - j(\bar{y}, z) \rangle \\
& \ + \psi_\beta(\bar{y} + j(\bar{y}, z)) - \psi_\beta(\bar{y}, z),
\end{split}
\end{equation}
and replacing this inequality into the definition of $\B^{\delta,1}$ we get
\begin{eqnarray*}
\B^{\delta, 1} 
&\leq& \epsilon^{-2} \int \limits_{B \setminus B_\delta} |j(\bar{x}, z) - j(\bar{y}, z)|^2 \nu(dz) 
+  \I^j[B\cap B_\delta^c](\psi_\beta, \bar{y})\\
&\leq&
C^2(1) \frac{|\bar{x}-\bar{y}|^2}{\epsilon^2}
+ \Lambda o_\beta(1)
\end{eqnarray*}
where we have used (MJ) with $R=1$ and Lemma~\ref{lemapsibeta} to control the nonlocal term applied to $\psi_\beta$. 

Thus, by~\eqref{comparison1}, we conclude that
\begin{equation}\label{Bdelta1comparison}
\B^{\delta, 1} \leq m_\beta (\epsilon) + o_\beta(1). 
\end{equation}

Now we address the estimate of $\B^{\delta,2}$. Note that
\begin{equation*}
\B^{\delta, 2} = \int_{B^c} [\bar{u}(\bar{x} + j(\bar{x}, z), \bar{t}) - v(\bar{y} + j(\bar{y}, z), \bar{t})
- (\bar{u}(\bar{x}, \bar{t}) - v(\bar{y}, \bar{t}))] \nu(dz).
\end{equation*}

For any $R>1,$ we divide this integral in two parts, 
$\B^{\delta, 2}= \B^{\delta, 2}[B_R^c]+\B^{\delta, 2}[B_R\setminus B],$ where 
the first one is integrated over $B_R^c$
and the second one over $B_R\setminus B.$

For the first one, we have
\begin{eqnarray*}
&& \B^{\delta, 2}[B_R^c]
\leq 
\Psi \int_{B_R^c}\nu(dz), 
\end{eqnarray*}
where we recall that $\Psi= 2 (||u||_{L^\infty(Q_T)} + ||v||_{L^\infty(Q_T)}).$

From (M), we can choose $R=R_\eta$ large enough in order that
the above term is less than $\eta/4.$

For the second term, using~\eqref{ineg-max127} we can get

\begin{eqnarray*}
\B^{\delta, 2}[B_{R_\eta}\setminus B]
&\leq &
\int \limits_{B_{R_\eta}\setminus B} 
[ \epsilon^{-2} |j(\bar{x}, z)-j(\bar{y}, z)|^2 
+2\epsilon^{-2}\langle\bar{x}-\bar{y}, j(\bar{x}, z)-j(\bar{y}, z)\rangle] \nu(dz)\\
&& \hspace*{1cm}
+  \int \limits_{B_{R_\eta}\setminus B} [\psi_\beta (\bar{y} + j(\bar{y}, z))-\psi_\beta (\bar{y})]\nu(dz)\\
&\leq& 2\epsilon^{-2} \int \limits_{B_{R_\eta\setminus B}} \Big{(}|j(\bar{x}, z)-j(\bar{y}, z)|^2
+ |\bar x - \bar y||j(\bar{x},z) - j(\bar{y},z)| \Big{)} \nu(dz) \\
&& \hspace*{1cm}
+ ||D\psi_\beta||_\infty  \int \limits_{B_{R_\eta}\setminus B} |j(\bar{y}, z)|\nu(dz)\\
&\leq&
(C^2(R_\eta)+C^1(R_\eta))  \epsilon^{-2} |\bar{x}-\bar{y}|^2+ \Lambda \beta C^0(R_\eta),
\end{eqnarray*}
by using (MJ) with $R={R_\eta}.$
Finally, we obtain
\begin{equation*}
\B^{\delta, 2} \leq \eta/4 + m_\beta(\epsilon) + o_\beta(1).
\end{equation*}

Using this and estimate~\eqref{Bdelta1comparison}, we get
\begin{equation}\label{B^deltacomparison}
\B^\delta \leq  \eta/4 +m_\beta(\epsilon) + o_\beta(1),
\end{equation}
where $m_\beta(\epsilon), o_\beta(1)$ depend on $\mu ,\eta$
but are independent of $\delta$.

\medskip

\noindent
5.- \textsl{Conclusion:} Joining~\eqref{Acomparison},~\eqref{Hcomparison},~\eqref{B_deltacomparison} and~\eqref{B^deltacomparison}, and using them 
into~\eqref{testingcomparison}, we conclude that
\begin{equation*}
-\epsilon^{-2}o_\delta(1) -m_\beta(\epsilon) -o_\beta(1) - \eta/2 \leq -\eta + o_\rho(1).
\end{equation*}

Finally, letting $\rho, \delta \to 0$ first, then $\epsilon \to 0$ for a fixed, small enough  $\beta$,
we obtain a contradiction since $\eta > 0.$ It follows that $M\leq 0$ and the proof is complete.
\qed

In what follows we discuss some important consequences of the comparison principle.

\subsection{Well-posedness.}
The consideration of the following boundedness condition 

\medskip

\noindent
{\bf (H0)} \textsl{There exists a constant $H_0 > 0$ such that $||H(\cdot, 0)||_\infty \leq H_0$}

\medskip

\noindent
allows us to provide existence to~\eqref{pareqIto}-\eqref{initialdata} via Perron's method.


\begin{cor}\label{corparcomparison}
Let $u_0 \in \mathrm{C}_b(\R^d)$, $A$ satisfying (A), $\I^j$ defined as in~\eqref{operator} in such a way its components satisfy assumptions (M), (MJ), and $H$ satisfying assumptions (H0), (H1), (H2). Then, there exists a unique viscosity solution 
$u \in C(\bar{Q})$ to problem~\eqref{pareqIto}-\eqref{initialdata}, which is also in $L^\infty(\bar{Q}_T)$ for all $T > 0$.
\end{cor}

The proof of this result follows classical arguments, see~\cite{Ishii3, usersguide}. It is posible to argument in $Q_T$ first and then extend it to the infinite time horizon. The role of the global sub and supersolution present in Perron's method is played by functions with the form $(x,t) \mapsto C_1t + C_2$, for suitable constants $C_1, C_2$ depending on the data and $T$. Uniqueness comes from Proposition~\ref{parcomparison}.

For the stationary case we can follow closely the previous arguments and obtain the analogous well-posedness result provided the equations is strictly proper.


\begin{prop}\label{comparisonstationary}
Let $\lambda > 0$, $A$ satisfying (A), $\I^j$ defined as in~\eqref{operator} in such a way its components satisfy assumptions (M), (MJ), and $H$ satisfying assumptions (H0), (H1), (H2).
Let $u$ be an USC bounded viscosity subsolution
and $v$ be a LSC bounded viscosity supersolution to equation~\eqref{eq}. Then, $u \leq v$ in $\R^d$. 

Moreover, if in addition we assume (H0), then there exists a unique viscosity solution $u \in C_b(\R^d)$ 
to equation~\eqref{eq}, for which we have the following bound
\begin{eqnarray}\label{uboundstationary}
||u||_\infty \leq \lambda^{-1} H_0.
\end{eqnarray}
\end{prop}


\subsection{Continuity results coming from comparison.} 
Comparison principle given by Proposition~\ref{parcomparison} allows us to obtain important continuity results for the solution of the addressed problems. The first result gives Lipschitz regularity in time if the initial data is smooth.


\begin{prop}\label{Lipschitztime} 
Consider the hypotheses of Proposition~\ref{parcomparison} and assume further 
that $u_0 \in C^2_b(\R^d)$ with $||u_0||_{C^2(\R^d)} < +\infty.$
Then, there exists a constant $\Lambda_0$ depending on the datas and $||u_0||_{C^2(\R^d)}$
such that each viscosity solution $u$ to problem~\eqref{pareqIto}-\eqref{initialdata} in 
$C(\bar{Q}) \cap L^\infty(\bar Q_T)$ for each $T > 0$ satisfies

\begin{equation*}
|u(x,t) - u(x, s)| \leq \Lambda_0 |t - s| \quad \mbox{for all} \ 0\leq s, t \leq T, \ x \in \R^d.
\end{equation*}
\end{prop}

\medskip

The proof of this result relies in comparison with functions with the form $(x, t) \mapsto \Lambda_0 t + u_0(x)$ for $\Lambda_0 \in \R$ adequate, together with translation invariance in time, see~\cite{Tchamba}.

A less direct consequence is stated in the following

\begin{prop}\label{propunifcont}
Let $u_0 \in BUC(\R^d)$ and assume that $\I^j$ and $H$ satisfy the assumptions of Corollary~\ref{corparcomparison}. Denote by $u \in C(\bar Q)$ the unique viscosity solution to~\eqref{pareqIto}-\eqref{initialdata} given in Corollary~\ref{corparcomparison}. Assume further that

\medskip
\noindent
(i) There exists a modulus of continuity $m_j$ such that for each $R > 0$ there exists $C_R > 0$ satisfying
\begin{equation*}
|j(x, z) - j(y, z)| \leq C_R \ m_j(|x - y|), \quad \mbox{for all} \ x, y \in \R^d, \ |z| \leq R. 
\end{equation*}

\medskip
\noindent
(ii) For each $R > 0$ there exists a constant $C_R > 0$ such that
\begin{equation*}
|H(x, p)| \leq C_R \quad \mbox{for all } \ x \in \R^d, \ |p| \leq R.
\end{equation*}

Then, there exists a modulus of continuity $m_T$ depending on the data, $T$ and $||u||_{L^\infty(Q_T)}$ such that
\begin{eqnarray*}
|u(x,t)-u(y,t)| \leq m_T(|x-y|) \quad \mbox{for all} \ x, y \in \R^d.
\end{eqnarray*}
\end{prop}

\medskip

\noindent
{\bf \textit{Proof:}} We only sketch the proof since most of the arguments are tedious but easily checkable.

By contradiction we assume that $u$ is not uniformly continuous in $x$. Then, there exist sequences $x_k, y_k\in\R^d,$
$t_k\in [0,T]$ such that $x_k-y_k$ tends to 0 and $u(x_k,t_k)-u(y_k,t_k)> \eta >0$.

We consider the function $v_k(x,t) = u(x+x_k,t)$. Denoting $u_{0, k}(x) = u_0(x + x_k)$, $A_k(x) = A(x + x_k)$, $j_k(x, z) = j(x + x_k, z)$ and $H_k(x, p) = H(x + x_k, p)$ we see that $v_k$ is a solution to the problem 
\begin{equation*}
\partial_t v_k - \I^{j_k}(v_k(\cdot, t), x) - \Tr(A_k D^2 v_k) + H_k(x, Dv_k) = 0 \quad \mbox{in} \ Q_T,
\end{equation*}
with initial data $v_k(\cdot, 0) = u_{0, k}$ in $\R^d$. 

The assumptions of the problem considered in the existence and uniqueness result together with the extra assumptions $(i)$ and $(ii)$ imply that the sequences of functions $\{ u_{k, 0} \}_k, \{ A_k \}_k, \{ H_k(\cdot, p) \}_k$, and $\{ j_k(\cdot, z) \}_k$ are locally uniformly continuous in $\R^d$ and bounded, for each $z, p \in \R^d$. Then, up to subsequences there exist functions $\tilde u_0, \tilde A, \tilde H$ and $\tilde j$ which are respective limit functions to the previous sequences, and the convergence is locally uniform in $\R^d$.
%
%
%
%
%
Then, defining the functions
\begin{equation*}
\bar v(x, t) = \limsup \limits_{y \to x, s \to t, k \to \infty} v_k(y, s), \quad 
\underline v(x, t) = \liminf \limits_{y \to x, s \to t, k \to \infty} v_k(y, s)
\end{equation*}
we apply the half-relaxed limits method (see~\cite{BP1}) to conclude that $\bar v, \underline v$ are, respectively, bounded viscosity sub and supersolutions to the problem
\begin{equation*}
\partial_t v - \I^{\tilde j}(v(\cdot, t), x) - \Tr(\tilde A D^2 v) + \tilde H(x, Dv) = 0 \quad \mbox{in} \ Q_T,
\end{equation*}
with the initial condition $v(\cdot, 0) = \tilde u_0$ in $\R^d$. Noticing that the limit problem satisfies the hypotheses of Proposition~\ref{parcomparison}, we conclude that $\bar v \leq \underline v$ and therefore they are equal.
%
But $u(x_k,t_k)-u(y_k,t_k)> \eta >0$ can be interpreted as 
$v_k(0,t_k)-v_k(y_k-x_k,t_k)> \eta >0$ and 
this would lead to 
a contradiction after taking limit as $k \to \infty$. This concludes the result.
\qed 

\medskip

Similar arguments can be given for the stationary problem~\eqref{eq} to conclude the following
\begin{prop}\label{propunifcont2}
Let $\lambda > 0$ and assume $\I^j, H$ satisfy 	the assumptions of Proposition~\ref{comparisonstationary}, and denote $u \in C(\R^d)$ the unique viscosity solution to~\eqref{eq}. Assume further conditions $(i), (ii)$ in Proposition~\ref{propunifcont}. Then, there exists a modulus of continuity $m$ depending on the data and $\lambda$ such that
\begin{eqnarray*}
|u(x)-u(y)| \leq m(|x-y|) \quad \mbox{for all} \ x, y \in \R^d.
\end{eqnarray*}
\end{prop}




\section{Lipschitz regularity - Bernstein method.}
\label{Lipschitzsection}

In this section we provide the Lipschitz regularity for solutions of problem~\eqref{eq}. This is accomplished by the introduction of a change 
of variables in the direction of Bernstein method, see~\cite{Barles2}.  As it can be seen in the literature, 
it is natural to assume Lipschitz regularity on the data to apply the mentioned method. Keeping assumptions 
(H1) and (A) as in the previous section, we require to strengthen assumptions (H2), (M) and (MJ) in the following sense.

\medskip

\noindent
{\bf (H2')} \textsl{Let $m$ as in (H1). There exists $L_H > 0$ and a modulus of continuity $\zeta$ 
such that, for all $x, y, p, q \in \R^d,$ $|q|\leq 1,$ we have
\begin{eqnarray*}
H(y, p + q) - H(x, p) \leq L_H |x - y| (1 + |p|^m) + \zeta(|q|)(1 + |p|^{m - 1}).
\end{eqnarray*}
}

\medskip

\noindent{\bf (M')} \textsl{ There exists $C_\nu > 0$ such that
\begin{equation*}
\int \limits_{\R^d} 1 \wedge |z|^2 \nu(dz) \leq C_{\nu},
\end{equation*}}

\noindent
{\bf (J1)} \textsl{ There exists a constant $C_j > 0$ such that for all $x, y, z \in \R^d$
\begin{equation*}
|j(x, z)| \leq C_j |z| \quad \mbox{and} \quad |j(x,z) - j(y,z)| \leq C_j |z| |x - y|.
\end{equation*}
}

\medskip

\noindent
{\bf (J2)} \textsl{ For each $a > 0$, there exists a constant $C_a > 0$ such that for all $x, y \in \R^d$
\begin{equation*}
\begin{split}
\int_{B_a^c} |j(x, z) - j(y, z)|\nu(dz) \leq C_a |x - y|
\end{split}
\end{equation*}
}

\begin{remark}\label{rmkLip}
Once we consider assumptions (A), (H1), together with (H2'), (M'), (J1) and (J2), we can apply all the results of the previous section since (M'), (J1) and (J2) imply (M),(MJ).
\end{remark}

The above assumptions are natural since they reflect the Lipschitz regularity of the data. 
This is evident in (J1) and (J2), meanwhile because of (H2') we have the function $x \mapsto H(x,0)$ is Lipschitz continuous.

Notice that (J2) reflects a compatibility condition among the measure $\nu$ and the jumps $j$. For example, (J2) is fullfilled by measures $\nu$ satisfying (M') and jumps satisfying
\begin{equation*}
|j(x,z) - j(y,z)| \leq C |x - y| \eta(z), \quad \mbox{for all } x, y \in \R^d, z \in B^c,
\end{equation*}
with $\eta(z)\nu(dz)$ finite in $B^c$.

\medskip

The main result of this paper is the following


\begin{teo}\label{teo1}{\bf (Lipschitz regularity)}
Let $\lambda \geq 0$. Assume $A$ satisfies (A), $H$ satisfies (H1), (H2'), and $\I^j$ defined in~\eqref{operator} such that $\nu, j$  satisfy assumptions (M'), (J1) and (J2). Let $u$ be a bounded uniformly continuous viscosity solution to problem~\eqref{eq}. 
Then, there exists $L > 0$ large enough such that
\begin{equation*}
|u(x) - u(y)| \leq L|x - y| \quad \mbox{for all} \ x,y \in \R^d.
\end{equation*}

The constant $L$ depends only on the data and $\mathrm{osc}(u)$.
\end{teo}


As we will see next, Theorem~\ref{teo1} is accomplished by an exponential change of variables inspired by the Bernstein method for viscosity solutions. In fact, we 
note that if $u \in C_b(\R^d)$ is a solution of~\eqref{eq}, replacing $u$ by $u - \inf \{u\} + 1$ we can assume $u \geq 1$ and then, under the change of variables $u = e^v$, we can prove that $v \geq 0$ satisfies the equation \footnote{In fact, $H$ should be changed into $H + \lambda(\inf \{u\} - 1)$. This new Hamiltonian has the same properties of $H$ since $\lambda \inf \{ u \}$ is bounded independently of $\lambda$ by~\eqref{uboundstationary}. Hence we keep the notation $H$ for the sake of simplicity.} 
\begin{equation}\label{expeq}
\lambda - \Tr(A D^2v) - \J^j(v, x) + e^{-v} H(x, e^v Dv) - |\sigma^T Dv|^2 = 0 \quad \mbox{in} \ \R^d,
\end{equation}
where $\J^j(v, x)$ is defined as
\begin{equation}\label{defJ}
\J^j(v,x) = \int \limits_{\R^d} [e^{v(x + j(x,z)) - v(x)} - 1 - \mathbf{1}_B(z) \langle Dv(x) , j(x,z) \rangle] \nu(dz).
\end{equation}

For simplicity, for $x, p \in \R^d$, $r \in \R$ we introduce the notation
\begin{equation}\label{defHtilde}
\tilde{H}(x, r, p) = \lambda  + e^{-r} H(x, e^r p) - |\sigma^T(x) p|^2.
\end{equation}

With this, using equation~\eqref{expeq} and the above notation, $u$ is a viscosity solution to~\eqref{eq} if and only if $v$ defined as $u = e^v$ is a viscosity solution to the equation
\begin{equation}\label{expeq2}
\tilde{H}(x, v(x), Dv(x)) - \Tr(A(x) D^2v(x)) - \J^j(v,x) = 0, \quad  x \in \R^d.
\end{equation}

By the above discussion, it is going to be convenient to argue over the equivalent equation~\eqref{expeq2}. As we mentioned in the introduction, of main interest is the treatment of the nonlocal term, and the nonlinearity arising in $\J^j$ after the exponential change of variables is an extra difficulty.

Now we present the proof of Theorem~\ref{teo1}. In its proof we use several technical estimates which are precisely stated and proved in the Appendix. 

\medskip

\noindent
{\bf \textit{Proof of Theorem~\ref{teo1}:}} Note that under the change $u = e^v$ 
with $u\geq 1,$ we have that $\mathrm{osc}(v)\leq \mathrm{osc}(u),$ and by the boundedness of $v$ we see that $v$ is still uniformly continuous. Then, if we get $v$ is Lipschitz continuous then we get the result for $u$. 

We will argue over the equation~\eqref{expeq2}. Assume by contradiction that for each $L > 0$ large enough (we can assume $L > 1$), there exists $\epsilon_L > 0$ such that
\begin{equation*}
\sup \limits_{\R^d \times \R^d} \{ v(x) - v(y) - L|x - y| \} \geq \ 2\epsilon_L.
\end{equation*}
and therefore, there exist $x_L, y_L \in \R^d$ such that
\begin{equation*}
v(x_L) - v(y_L) - L|x_L - y_L| \geq \epsilon_L. 
\end{equation*}

At this point we consider a nonnegative function $\psi \in C^2_b(\R^d)$ satisfying assumptions~\eqref{defpsi} with $\Psi = \mathrm{osc} (v)$ and $\Lambda > 0$ depending only on $\mathrm{osc}(v)$, and for $\beta > 0$ we consider $\psi_\beta$ as in~\eqref{psibeta}. Then, setting $\phi(x,y) := L|x - y| + \psi_\beta(y)$ we have
\begin{equation}\label{suplocalized}
\mathop{\rm sup}_{\R^d \times \R^d} \{ v(x) - v(y) - \phi(x,y) \}
\geq 
v(x_L) - v(y_L) - L|x_L - y_L| - \psi_\beta(y_L) \geq \epsilon_L >0
\end{equation}
for all $\beta$ small enough to have $1/\beta >|y_L|.$ Moreover, 
since $\psi_\beta = \mathrm{osc}(v)$ in $B_{2/\beta}^c$,
the supremum in~\eqref{suplocalized} is achieved 
at some point $(\bar{x}, \bar{y}) \in \R^d \times \R^d$, 
with $\bar{x} \neq \bar{y}$ for each $L, \beta > 0$, and we get
\begin{equation}\label{x-yabove}
L|\bar x - \bar y| \leq \osc(v). 
\end{equation}

Using Proposition~\ref{propunifcont}, let $\omega: [0,+\infty) \to [0,+\infty)$, $\omega(0) = 0$, be an increasing modulus of continuity for $v$. Then we can write
\begin{equation*}
|v(x) - v(y)| \leq \omega(|x - y|) \quad \mbox{for all}\ x, y \in \R^d.
\end{equation*}

By~\eqref{suplocalized} we have $\epsilon_L \leq v(\bar{x}) - v(\bar{y})$ and therefore
\begin{equation}\label{x-ybelow}
0 < \omega^{-1}(\epsilon_L) \leq |\bar{x} - \bar{y}|, 
\end{equation}
which is a lower bound for $|\bar{x} - \bar{y}|$, uniform in terms of $\beta$.

Now we would like to use the viscosity inequality for $v$ at $\bar{x}$ and $\bar{y}$. 
Using Ishii-Jensen lemma for nonlocal equations, for all $\delta > 0$ and all $\rho > 0$
small enough we have
\begin{equation}\label{testingLip}
\begin{split}
\tilde{H}(\bar{x}, v(\bar{x}), L\bar{p})  - \Tr(A(\bar{x}) X_\rho)  & \\
- \J^j[B_\delta](\phi(\cdot, \bar{y}), \bar{x}) 
- \J^j[B_\delta^c](v, \bar{x}, L\bar{p}) & \leq o_\rho(1), \\
\tilde{H}(\bar{y}, v(\bar{y}), L\bar{p} - \bar{q})  - \Tr(A(\bar{y}) Y_\rho)  & \\
- \J^j[B_\delta](-\phi(\bar{x}, \cdot), \bar{y}) 
- \J^j[B_\delta^c](v, \bar{y}, L\bar{p} - \bar{q}) & \geq -o_\rho(1), 
\end{split}
\end{equation}
where $\bar{p} = (\bar{x} - \bar{y})/|\bar{x} - \bar{y}|$ and $\bar{q} = D\psi_\beta(\bar{y})$. The matrices 
$X_\rho, Y_\rho \in \mathbb{S}^{d}$ satisfy the inequality
\begin{equation}\label{IJineq}
-\rho^{-1} I_{2d} \leq \left [ \begin{array}{cc} X_\rho & 0 \\ 0 & -Y_\rho \end{array} \right ] \leq D^2\phi(x,y) + o_\rho(1). 
\end{equation}

We remark that all the terms $o_\rho(1)$ arising in~\eqref{testingLip} and~\eqref{IJineq} satisfy $o_\rho(1) \to 0$ if $L, \beta > 0$ are fixed.

Subtracting both inequalities in~\eqref{testingLip}, we get
\begin{equation}\label{HBALip}
\mathcal{H} \leq \mathcal{A} + \B_\delta + \B^\delta, 
\end{equation}
where
\begin{equation}\label{termes123}
\begin{split}
\mathcal{H} = & \ \tilde{H}(\bar{x}, v(\bar{x}), L\bar{p}) - \tilde{H}(\bar{y}, v(\bar{y}), L\bar{p} - \bar{q}), \\
\B_\delta = & \ \J^j[B_\delta](\phi(\cdot, \bar{y}), \bar{x}) 
- \J^j[B_\delta](-\phi(\bar{x}, \cdot), \bar{y}), \\
\B^\delta = & \ \J^j[B_\delta^c](v, \bar{x}, L\bar{p}) 
- \J^j[B_\delta^c](v, \bar{y}, L\bar{p} - \bar{q}), \\
\mathcal{A} = & \ \Tr(A(\bar{x}) X_\rho) - \Tr(A(\bar{y}) Y_\rho) + o_\rho(1). 
\end{split}
\end{equation}

In what follows, we estimate each term present in~\eqref{HBALip} in order to get the desired contradiction 
by taking $L$ large enough in terms of the data.


\medskip

\noindent
1.- \textsl{Estimates for $\mathcal{A}$.} Note that in this case we have
\begin{equation*}
\phi(x,y) = L \varphi_1(x,y) + \psi_\beta(y),
\end{equation*}
where $\varphi_1$ is defined in~\eqref{varphi}.

Hence, using the estimate given by Lemma~\ref{lemaADvarphi} and applying the estimates for the second derivatives of $\psi_\beta$
in Lemma~\ref{lemapsibeta}, we conclude that
\begin{equation*}
\begin{split}
\Tr(A(\bar{x}) X_\rho - A(\bar{y}) Y_\rho) 
\leq  2 L L_\sigma^2 |\bar{x} - \bar{y}| + \beta^2 L_\sigma^2 \Lambda + o_\rho(1)
\end{split}
\end{equation*}
and therefore
\begin{equation}\label{ALip}
\begin{split}
\mathcal{A} \leq 2 L_\sigma^2 L|\bar{x} - \bar{y}| + o_\beta(1) + o_\rho(1).
\end{split}
\end{equation}


\medskip

\noindent
2.- \textsl{Estimates for $\mathcal{H}$.} By definition of $\tilde{H}$ in $\mathcal{H}$ we have
\begin{equation*}\label{AteoLip1}
\mathcal{H} = \mathcal{H}_1  + \mathcal{H}_2,
\end{equation*}
with
\begin{equation*}
\begin{split}
\mathcal{H}_1 = & \ e^{-v(\bar{x})} H(\bar{x}, e^{v(\bar{x})}L\bar{p}) - e^{-v(\bar{y})} H(\bar{y}, e^{v(\bar{y})}(L\bar{p} - \bar{q})) \\
\mathcal{H}_2 = & \ -|\sigma^T(\bar{x})L\bar{p}|^2 + |\sigma^T(\bar{y})(L\bar{p} - \bar{q})|^2.
\end{split}
\end{equation*}

For $\mathcal{H}_2$, using the estimates for the gradient of $\psi_\beta$ and the fact that $|\bar{p}| = 1$, we can write
\begin{equation*}
\begin{split}
\mathcal{H}_2 = & \ L^2 ( |\sigma^T(\bar{y})(\bar{p} - L^{-1} \bar{q})|^2 - |\sigma^T(\bar{x})\bar{p}|^2) \\ 
\geq & \ L^2 ( |\sigma^T(\bar{y})\bar{p}|^2 - |\sigma^T(\bar{x})\bar{p}|^2 
- 2L^{-1} \beta^2 L_\sigma^2  \Lambda ),
\end{split}
\end{equation*}
and applying the boundedness and the Lipschitz estimates of $\sigma$ given in (A), we conclude
\begin{equation}\label{H2}
\mathcal{H}_2 \geq -2L_\sigma^2 L^2  |\bar{x} - \bar{y}| 
- 2L_\sigma^2\Lambda^2\tilde{C} \beta^2 L.
\end{equation}

The estimate from below of $\mathcal{H}_1$ is a crucial step in the weak
Bernstein method. It is based on the superlinearity of the Hamiltonian which is encoded
in Assumption (H1) (see~\eqref{Hcoercive}). We
define $\mu = e^{v(\bar{y}) - v(\bar{x})}$ and we notice that 
$0 < \mu < 1$ since $v(\bar y) < v(\bar x)$ for all $L, \beta$. 
Using this we can write
\begin{equation*}
\mathcal{H}_1 = e^{-v(\bar{y})} \Big{(} \mu H(\bar{x}, \mu^{-1} e^{v(\bar{y})} L\bar{p})
- H(\bar{y}, e^{v(\bar{y})} (L \bar{p} - \bar{q})) \Big{)}.
\end{equation*}

Applying (H1), (H2') with $\beta$ small enough in order that
$|e^{v(\bar{y})}\bar{q}|\leq 1$ and recalling that $m > 1$ and $v \geq 0$, we can write
\begin{equation*}
\begin{split}
\mathcal{H}_1 
\geq & \ e^{-v(\bar{y})} \Big{(} (1 - \mu) (b_m L^m e^{mv(\bar{y})} - K) - L_H |\bar{x} - \bar{y}|(1 + L^m e^{m v (\bar{y})}) \\
 & \qquad \quad - \zeta(|e^{v(\bar{y})}\bar{q}|)
L^{m - 1} e^{(m - 1)v(\bar{y})}  -\zeta(|e^{v(\bar{y})}\bar{q}| \Big{)} \\
\geq & \ L^m e^{(m - 1)v(\bar{y})} \Big{(} (b_m - KL^{-m})(1-\mu) - L_H |\bar{x} - \bar{y}|(1 + L^{-m}) \Big{)} - o_\beta(1),
\end{split}
\end{equation*}
where $o_\beta(1) \to 0$ as $\beta \to 0$, but depending on $L$ and $||v||_\infty$. 
Taking $L$ satisfying
\begin{eqnarray}\label{choixL1}
L^{m} \geq \max \{1, 2^{-1} b_m^{-1} K\},
\end{eqnarray}
we conclude that
\begin{equation*}
\mathcal{H}_1 \geq L^m e^{(m - 1)v(\bar{y})} \Big{(} 
\frac{b_m}{2}(1-\mu) - 2L_H |\bar{x} - \bar{y}| \Big{)} - o_\beta(1).
\end{equation*}

Now, we have
\begin{eqnarray*}
1-\mu = e^{v(\bar{y}) - v(\bar{x})}( e^{v(\bar{x}) - v(\bar{y})}-1)
\geq e^{v(\bar{y}) - v(\bar{x})}(v(\bar{x}) - v(\bar{y}))
\geq e^{-\mathrm{osc}(v)}L|\bar{x} - \bar{y}|
\end{eqnarray*}
since $\mathrm{osc}(v)\geq v(\bar{x}) - v(\bar{y})\geq L|\bar{x} - \bar{y}|$
by maximality of $(\bar{x},\bar{y})$ in~\eqref{suplocalized}. 
From this, taking 
\begin{eqnarray}\label{choixL2}
L \geq 8 L_H e^{\mathrm{osc}(v)}b_m^{-1}
\end{eqnarray}
and using that $v \geq 0$ and $m > 1$, we conclude that
\begin{equation*}
\mathcal{H}_1 \geq 
\frac{b_m}{4} e^{-\mathrm{osc}(v)} L^{m + 1} |\bar{x} - \bar{y}| - o_\beta(1).
\end{equation*}

Recalling that $\mathcal{H} = \mathcal{H}_1 + \mathcal{H}_2$, we join the last estimate 
and~\eqref{H2} to obtain
\begin{equation*}
\mathcal{H} \geq  (\frac{b_m}{4} e^{-\mathrm{osc}(v)}
L^{m + 1} - 2L_\sigma^2 L^2) |\bar{x} - \bar{y}| - o_\beta(1),
\end{equation*}
and therefore, since $m > 1$, taking 
\begin{eqnarray}\label{choixL3}
L^{m-1}\geq 16L_\sigma^2 b_m^{-1} e^{\mathrm{osc}(v)},
\end{eqnarray}
we finally conclude that
\begin{equation}\label{HLip}
\mathcal{H} \geq c L^{m + 1} |\bar{x} - \bar{y}| - o_\beta(1),
\end{equation}
where $c = b_m e^{-\mathrm{osc}(v)}/8 > 0$ 
and $o_\beta(1) \to 0$ as $\beta \to 0$ for fixed $L$.


\medskip

\noindent
3.-\textsl{Estimates for $\B_\delta$.} 
At this point, recalling that  $|\bar{x} - \bar{y}|>0$ by~\eqref{x-ybelow},
we are going to consider $\delta > 0$ small enough depending on $|\bar{x} - \bar{y}|$
such that $|j(\bar{x}, z)|, |j(\bar{y}, z)| \leq |\bar x -\bar y |/2$ for each $|z| \leq \delta$.
This is possible by assumption~(J1).

Thus, we have $x \mapsto |\bar x - \bar y + x|$, $y \mapsto |\bar x - \bar y - y|$ are smooth in $B_\delta$.
Then, applying Lemma~\ref{lemaIexp} we can write
\begin{equation*}
\begin{split}
\B_\delta \leq \I^j[B_\delta](\phi(\cdot, \bar{y}), \bar{x}) 
+ \I^j[B_\delta](\phi(\bar{x}, \cdot), \bar{y}) + C(L^2+ \beta^2 \Lambda^2) o_\delta(1),
\end{split}
\end{equation*}
and by definition of $\phi$ together with Lemma~\ref{lemapsibeta} and assumptions (M'), (J1) we arrive at
\begin{equation}\label{BdeltaLip}
\B_\delta \leq \Big{(} L|\bar{x} - \bar{y}|^{-1} + L^2 + o_\beta(1) \Big{)} o_\delta(1).
\end{equation}

Note that the estimate above is appropriate since we will send
$\delta\to 0$ first in the global estimate in Step 5.
\medskip

\noindent
4.- \textsl{Estimates for $\B^\delta$.} 
This is the key new estimate to perform the weak Bernstein's method in the nonlocal case and therefore we state it as a lemma.
\begin{lema}\label{estim-nonlocale}
There exists $C > 0$ depending only on the datas $C_j, C_{\nu, j}, C_\nu$ and on $\mathrm{osc}(v)$ such that
\begin{equation}\label{B^deltaLip}
\mathcal B^\delta \leq C L^2|\bar x- \bar y| + o_\beta(1),
\end{equation}
where $o_\beta(1) \to 0$ as $\beta \to 0$ for $L > 0$ fixed.
\end{lema}

Then, the rest of Step 4 is devoted to the proof of this lemma.
We start with the following notation: for a function $f: \R^d \to \R$ and $x, z \in \R^d$, we denote
\begin{equation*}
\Delta_x f = \Delta_xf(z) = f(x + j(x,z)) - f(x).
\end{equation*}

We also consider $\Theta_i(z)$ for $i = 1,2,3$, $z \in \R^d$, defined as
\begin{equation*}
\begin{split}
& \Theta_1(z) = L(|\bar{x} - \bar{y} + j(\bar{x},z)| - |\bar{x} - \bar{y}|), \\ 
& \Theta_2(z) =  -L(|\bar{x} - \bar{y} - j(\bar{y},z)| - |\bar{x} - \bar{y}|) - \Delta_{\bar{y}}\psi_\beta, \\ 
& \Theta_3(z) = L(|\bar{x} - \bar{y} + j(\bar{x},z) - j(\bar{y},z)| - |\bar{x} - \bar{y}|) + \Delta_{\bar{y}}\psi_\beta, 
\end{split}
\end{equation*}
and notice that the maximality of $(\bar{x}, \bar{y})$ in~\eqref{suplocalized} 
implies, for each $z \in \R^d$, the following inequalities
\begin{equation}\label{ineqDelta}
\begin{split}
\Delta_{\bar{x}}v \leq \Theta_1, \quad
\Delta_{\bar{y}}v \geq  \Theta_2, \quad
\Delta_{\bar{x}}v - \Delta_{\bar{y}}v \leq \Theta_3.
\end{split}
\end{equation}

We write $\B^\delta=\B_1+\B_2^\delta$ with
\begin{equation*}
\begin{split}
\B_1 = & \ \int_{B^c} [e^{\Delta_{\bar{x}}v} - e^{\Delta_{\bar{y}}v}] \nu(dz) \\
\B^\delta_2 = & \ \int_{B \setminus B_{\delta}} [e^{\Delta_{\bar{x}}v} - e^{\Delta_{\bar{y}}v} 
- L \langle \bar{p}, j(\bar{x}, z) - j(\bar{y}, z) \rangle - \langle \bar{q}, j(\bar{y}, z) \rangle] \nu(dz).
\end{split}
\end{equation*}

Our aim is to estimate from above the integrals. Of main importance is to get estimates for the integral over $\B^\delta_2$ independent of $\delta$ since we are going to take $\delta \to 0$ first at the end of this proof.

We start with the estimate of $\B_1$. We first remark that we can integrate only on the set $\mathcal{P}_1$
where $e^{\Delta_{\bar{x}}v} - e^{\Delta_{\bar{y}}v} \geq 0$, i.e., 
where $\Delta_{\bar{x}}v - \Delta_{\bar{y}}v \geq 0$. 

On this set, we have
\begin{eqnarray*}
e^{\Delta_{\bar{x}}v} - e^{\Delta_{\bar{y}}v}= e^{\Delta_{\bar{x}}v}(1- e^{\Delta_{\bar{y}}v-\Delta_{\bar{x}}v})
\leq e^{\mathrm{osc}(v)}(\Delta_{\bar{x}}v-\Delta_{\bar{y}}v)
\leq e^{\mathrm{osc}(v)} \Theta_3
\end{eqnarray*}
since $\Delta_{\bar{x}}v\leq \mathrm{osc}(v)$ and $0\leq 1-e^{-r}\leq r$ for $r\geq 0.$

Noticing that $\Theta_3 \leq L| j(\bar{x},z) - j(\bar{y},z)| + \Delta_{\bar{y}}\psi_\beta$ 
and using (J2) and Lemma~\ref{lemapsibeta}, we get
\begin{eqnarray*}
\B_1 \leq
 e^{\mathrm{osc}(v)} \left(L \int_{B^c\cap \mathcal{P}_1} | j(\bar{x},z) - j(\bar{y},z)| \nu(dz)
+ \I^j [B^c\cap \mathcal{P}_1](\psi_\beta , \bar{y})\right),
\end{eqnarray*}
and using (J2) for the first integral and Lemma~\ref{lemapsibeta} for the second term in the righ-hand side, we arrive at
\begin{equation}\label{B1delta}
\B_1 \leq e^{\mathrm{osc}(v)} (L C_1 | \bar{x} - \bar{y} |  + o_\beta (1)),
\end{equation}
where $C_1$ is given by (J2) for $a=1$ and
we point out that $o_\beta(1) \to 0$ as $\beta \to 0$ uniformly in all the other variables.

Now we deal with the estimate of $\B_2^\delta$. A key fact is that it is enough to integrate 
$$
\Psi(z):= e^{\Delta_{\bar{x}}v} - e^{\Delta_{\bar{y}}v} - L \langle \bar{p}, j(\bar{x}, z) - j(\bar{y}, z) \rangle - \langle \bar{q}, j(\bar{y}, z) \rangle
$$
only on the set  $\mathcal{P}_2$ where $\Psi(z)\geq 0.$ This consideration allows us to get the relevant estimates to apply the ``linearization" procedure provided by Lemma~\ref{enlt1}, which is proven in the Appendix.
In fact, notice that by the third inequality in~\eqref{ineqDelta},
applying (J1) and the properties of $\psi_\beta$, for each $z$ we can write
\begin{equation*}
\Delta_{\bar x} v  - \Delta_{\bar y} v \leq C(L |\bar x - \bar y| + o_\beta(1))|z|,
\end{equation*}
where $C > 0$ depends only on the data. 

On the other hand, for each $z \in \mathcal P_2$ we can write
\begin{equation*}
-(LC_j | \bar{x} - \bar{y} | +o_\beta(1))|z| 
\leq L \langle \bar{p}, j(\bar{x}, z) - j(\bar{y}, z) \rangle + \langle \bar{q}, j(\bar{y}, z) \rangle 
\leq e^{\Delta_{\bar{x}}v}-  e^{\Delta_{\bar{y}}v},
\end{equation*}
where the first inequality comes from (J1) and the properties of $\psi_\beta$, and the second
one from the definition of $\mathcal P_2$. Now, since
\begin{equation*}
e^{\Delta_{\bar{x}}v}-  e^{\Delta_{\bar{y}}v}
= e^{\Delta_{\bar{x}}v}(1-e^{\Delta_{\bar{y}}v-\Delta_{\bar{x}}v})
\leq e^{\Delta_{\bar{x}}v}(\Delta_{\bar{x}}v-\Delta_{\bar{y}}v),
\end{equation*}
we join this inequality and the previous one to conclude that
\begin{equation*}
-Ce^{\osc (v)}(L| \bar{x} - \bar{y} | +o_\beta(1))|z| \leq \Delta_{\bar x}v - \Delta_{\bar y}v
\end{equation*}
where $C > 0$ depends only on the data. Thus, the upper and lower bound for 
$\Delta_{\bar x}v - \Delta_{\bar y}v$ can be summarized as
\begin{equation}\label{ineqDelta1}
|\Delta_{\bar x}v - \Delta_{\bar y}v| \leq C(L|\bar x - \bar y| + o_\beta(1))|z| \quad \text{for $z\in  \mathcal{P}_2$}
\end{equation}
where $C$ denotes a constant which may vary line to line
but depends only the data and $\mathrm{osc}(v)$.

Now, using the first and second inequalities in~\eqref{ineqDelta} together with (J1) and the properties of $\psi_\beta$ we can write
\begin{equation*}
\Delta_{\bar x} v \leq  LC_j |z|, \quad \Delta_{\bar y} v 
\geq - \left(LC_j + o_\beta(1)\right)|z|.
\end{equation*}

These inequalities and~\eqref{ineqDelta1} allows us to obtain, for $z \in \mathcal{P}_2$
\begin{equation*}
|\Delta_{\bar{x}}v |, |\Delta_{\bar{y}}v| \leq C\left( L + L|\bar x - \bar y| + o_\beta(1)\right)|z|.
\end{equation*}

Finally, since we can assume $L > 1$ and by~\eqref{x-yabove} we conclude that
\begin{equation}\label{*}
|\Delta_{\bar{x}}v |, |\Delta_{\bar{y}}v| \leq C\left(L + o_\beta(1)\right)|z|.
\end{equation}

Recalling that we also have that $|\Delta_{\bar{x}}v |, |\Delta_{\bar{y}}v| \leq \mathrm{osc} (v)$, in view of~\eqref{ineqDelta1} and~\eqref{*} we can apply Lemma~\ref{enlt1} with $g(x,z) = \Delta_x v(z)$, $C_1$ just depending on the data, $C_2 = \mathrm{osc}(v), b = o_\beta(1)$ and $\mathcal P = \mathcal P_2$ to conclude that
\begin{equation*} 
e^{\Delta_{\bar x}v} - e^{\Delta_{\bar y}v} \leq \Delta_{\bar x}v - \Delta_{\bar y}v
+ C (L^2 |\bar x - \bar y| + o_\beta(1))|z|^2, \quad \mbox{for} \ z \in \mathcal P_2.
\end{equation*}

%
%
%
%
%

Then, by using this estimate and the last inequality in~\eqref{ineqDelta}, for each $z \in \mathcal P_2$ we get
\begin{equation*}
\Psi(z) \leq \Psi_1(z) + \Psi_2(z),
\end{equation*}
where
\begin{equation*}
\begin{split}
\Psi_1(z) = &  \ L\Big{(} |\bar{x} - \bar{y} + j(\bar{x},z) - j(\bar{y},z)| - |\bar{x} - \bar{y}| - 
\langle \bar{p}, j(\bar{x}, z) - j(\bar{y}, z)\rangle \Big{)}, \\
\Psi_2(z) = & \ \Delta_{\bar{y}}\psi_\beta - \langle \bar{q}, j(\bar{y}, z)  \rangle 
+ C (L^2 |\bar x - \bar y| + o_\beta(1))|z|^2.
\end{split}
\end{equation*}

For the integral term concerning $\Psi_2,$ we note that the integral of the first two terms is exactly
$\I^j[(B\setminus B_\delta)\cap \mathcal{P}_2](\psi_\beta , \bar{y})$ which is $o_\beta(1)$
by Lemma~\ref{lemapsibeta}. Using (M') to estimate the last term, we infer the existence
of a constant $C > 0$ not depending on $\delta, \beta$ or $ L$ such that
\begin{equation*}
\int_{\mathcal P_2 \cap (B \setminus B_\delta)} \Psi_2(z)\nu(dz) \leq C (L^2 |\bar x - \bar y| + o_\beta(1)),
\end{equation*}
where $o_\beta(1) \to 0$ when $L > 0$ is fixed.

For the estimate of the integral term related to $\Psi_1$ we use the estimate given by Lemma~\ref{estim-nonlocale1} proven in the Appendix to conclude that
\begin{equation*}
\int_{\mathcal P_2 \cap (B \setminus B_{\delta})} \Psi_1(z) \nu(dz) \leq C L |\bar x - \bar y|,
\end{equation*}
for some $C > 0$ depending on the data.

In view of the above estimates and since we assume $L > 1$, we finally arrive at
\begin{equation*}
\B_2^\delta \leq CL ^2 |\bar{x} - \bar{y}| + o_\beta(1).
\end{equation*}

Putting together this last estimate and the estimate for $\B_1$ in~\eqref{B1delta} 
we finally obtain~\eqref{B^deltaLip} as stated in Lemma~\ref{estim-nonlocale}.

\medskip

\noindent
5.- \textsl{Conclusion of the proof of the Theorem.} 
Replacing~\eqref{ALip},~\eqref{HLip},~\eqref{BdeltaLip} and~\eqref{B^deltaLip} into~\eqref{HBALip}, we obtain
\begin{eqnarray}\label{finalBernstein}
&& c L^{m + 1} |\bar{x} - \bar{y}| \leq C L^2 |\bar{x} - \bar{y}|\\
\nonumber 
&& \hspace*{3cm} + ( L|\bar{x} - \bar{y}|^{-1} + L^2 +  o_\beta(1) ) o_\delta(1) 
+ o_\beta(1) + o_\rho(1),
\end{eqnarray}
where $o_\rho(1) \to 0$ as $\rho \to 0$, $o_\delta(1) \to 0$ as $\delta \to 0$ uniformly in the remaining variables, 
$o_\beta(1) \to 0$ as $\beta \to 0$ when $L > 0$ is fixed and 
$L, C, c > 0$ depend only on the datas $C_\nu, C_{\nu,j}, C_j, L_\sigma, L_H, b_m, K$
and $\mathrm{osc}(v)$. 

More precisely, we can fix $L$ at the beginning such that~\eqref{choixL1},~\eqref{choixL2},~\eqref{choixL3} hold and in addition
\begin{equation*}
L \geq (c^{-1} C)^{1/(m - 1)} + 1.
\end{equation*}

With this choice and since $|\bar{x} - \bar{y}|$ is uniformly positive in terms of $\beta$ by~\eqref{x-ybelow}, 
making $\rho\to 0, \delta \to 0$, $\beta \to 0$ we arrive at a contradiction with~\eqref{finalBernstein}, which ends the 
proof of the theorem.
\qed
\bigskip




We can extend the Lipschitz regularity 
in the $x$ variable for the solution to the parabolic problem~\eqref{pareqIto}-\eqref{initialdata}.


\begin{prop}\label{Lipspaceparabolic}
Let $u_0$ be bounded and Lipschitz function in $\R^d$ with Lipschitz constant $L_0 > 0$. Assume $A$, $H$ and $\I^j$ defined in~\eqref{operator} satisfy the assumptions of Theorem~\ref{teo1}. Let $u \in C(\bar Q)$ be the unique viscosity solutions  to to problem~\eqref{pareqIto}-\eqref{initialdata} given by Corollary~\ref{corparcomparison}.


Then, there exists $L > 0$ depending on the data, $L_0$ and 
\begin{equation*}
\mathrm{osc}_T (u):= \mathop{\rm sup}_{t\in [0,T]} \mathrm{osc}(u(\cdot,t))=
 \mathop{\rm sup}_{t\in [0,T]} \{\mathop{\rm sup}_{\R^d} u(\cdot,t) - \mathop{\rm inf}_{\R^d} u(\cdot,t)\}
\end{equation*}
such that
\begin{equation*}
|u(x, t) - u(y, t)| \leq L |x - y|, \quad \mbox{for all} \ x, y \in \R^d, \ t \in [0,T].
\end{equation*}
\end{prop}


\noindent
{\bf \textit{Sketch of the Proof:}} As in the proof of Theorem~\ref{teo1}, 
we may assume without loss of generality that $u\geq 1$ and
we argue over the function $v$ defined though the change of variables
$u(x,t) = e^{v(x,t)}$ for all $(x,t) \in \bar{Q}$. 
Hence, proving the Lipschitz continuity in $x$ for $v$, we 
conclude the  desired property for $u$. 

We start by proving the result for $u_0\in C^2(\R^d)$ with
$||u_0||_{C^2(\R^d)} < +\infty$  in order to be able
to use Proposition~\ref{Lipschitztime}.

The new function $v$ solves the problem
\begin{equation*} 
\left \{ \begin{array}{rll}
\partial_t v - \mathrm{Tr}(AD^2v) - \J^j(v(\cdot, t), x) + \tilde{H}(x, v, Dv) &=  0 \quad & \mbox{in} \ Q \\
v(\cdot, 0) &= e^{u_0} \quad & \mbox{in} \ \R^d,
\end{array} \right .
\end{equation*}
where $\J_x^j$ is defined in~\eqref{defJ} and $\tilde{H}$ is defined by~\eqref{defHtilde} with $\lambda=0.$

As in Theorem~\ref{teo1}, we argue by contradiction, assuming that for all $L \geq 1$ 
large enough, there exists $x_L,y_L\in\R^d,$ $t_L\in [0,T]$ and
$\epsilon_L > 0$ such that
\begin{equation*}
\begin{split}
& \sup \limits_{(x,y,t) \in \R^d \times \R^d \times [0,T]} \!\!\!\!\!\!\!\!\!\!\!
\{ v(x,t) - v(y, t) - L|x - y|\} \\
& \qquad \geq \ v(x_L,t_L) - v(y_L, t_L) - L|x_L - y_L|\geq 
2 \epsilon_L.
\end{split}
\end{equation*}

Hence, introducing the localization function $\psi_\beta$ defined in~\eqref{psibeta} 
as in Theorem~\ref{teo1} but replacing $\mathrm{osc} (v)$ by $\mathrm{osc}_T (v)$,
for all $\beta >0$ small enough with respect to $L$ and all $\eta > 0$ we have
\begin{equation}\label{parLipM}
\sup \limits_{(x,y,s,t) \in (\R^d)^2 \times [0,T]^2} \{ v(x,s) - v(y, t) - L|x - y| - \psi_\beta(y) - \eta^{-1}(s - t)^2 \} 
\geq 2 \epsilon_L.
\end{equation}

Applying Proposition~\ref{Lipschitztime} with a constant $\Lambda_0 > 0$ depending on $||u_0||_{C^2(\R^d)}$ and using the definition of $\osc_T$ we can write
\begin{equation*}
\begin{split}
& v(x,s) - v(y, t) - L|x - y| - \psi_\beta(y) - \eta^{-1}(s - t)^2 \\
\leq & \ \osc_T(v) - \psi_\beta(y) - L|x - y| + \Lambda_0|s-t| - \eta^{-1}(s - t)^2.
\end{split} 
\end{equation*}

Notice that $\Lambda_0 a - \eta^{-1} a^2 \leq \Lambda_0^2 \eta/4$ for all $a > 0$. Using this and the properties of $\psi_\beta$ we get the inequality
\begin{equation*}
v(x,s) - v(y, t) - L|x - y| - \psi_\beta(y) - \eta^{-1}(s - t)^2 \leq -L|x - y| + \Lambda_0^2 \eta /4
\end{equation*}
for all $|y|\geq 2/\beta$. It follows that for all $\eta$ small enough in terms of $\epsilon_L$ and $\Lambda_0$,
the supremum in~\eqref{parLipM} 
is achieved at some point $(\bar{x}, \bar{y}, \bar{s}, \bar{t}) \in (\R^d)^2 \times [0,T]^2$
and $\bar{x}, \bar{y}$ are in a bounded set depending only on $\mathrm{osc}_T (v)$
and $\beta.$ Moreover, by classical results, 
$\eta^{-1}(\bar{s} - \bar{t})^2 \to 0$ as $\eta \to 0$ uniformly in $L,\beta.$

In particular, up to take a subsequence $\eta\to 0,$ we get
$\bar{x}\to \bar{x}^*,$  $\bar{y}\to \bar{y}^*$ and $\bar{s}, \bar{t} \to t^*.$ 
Choosing $L > L_0$ implies $t^* > 0$ since otherwise passing to the limit as 
$\eta \to 0$ in~\eqref{parLipM} we arrive to a contradiction with $u_0$ $L_0$-Lipschitz
continuous. It follows that $\bar{s}, \bar{t}>0$ for sufficiently small $\eta.$
Likewise, we obtain that $\bar{x}\not= \bar{y}$ for $\eta$ sufficiently small.

It follows that, for $\eta$ small enough, we can write the viscosity inequalities
for the subsolution $v$ at $(\bar{x}, \bar{s})$ and  
for the supersolution $v$ at $(\bar{y}, \bar{t})$
to obtain exactly~\eqref{HBALip}, that is
\begin{equation}\label{nlleHBALip} 
\mathcal{H} \leq \mathcal{A} + \B_\delta + \B^\delta, 
\end{equation}
where $v(\bar{x})$ is replaced by $v(\bar{x}, \bar{s})$
and $v(\bar{y})$ by $v(\bar{y}, \bar{t})$ in~\eqref{termes123}.

Sending $\eta\to 0,$ we still have~\eqref{nlleHBALip}
where $\bar{x}$ is replaced by $\bar{x}^*,$  
$\bar{y}$ by $\bar{y}^*,$ $\bar{s}, \bar{t}$ by $t^*$
and $(\bar{x}^*,\bar{y}^*, t^*)$ is a maximum point of
\begin{equation}\label{parLipMsanseta}
\sup \limits_{(x,y,t) \in (\R^d)^2 \times [0,T]} \{ v(x,s) - v(y, t) - L|x - y| - \psi_\beta(y)\} 
\geq 2\epsilon_L.
\end{equation}

By compactness of $[0,T],$ we can assume that $t^*\to \hat{t}$ as
$\beta \to 0.$ Using Proposition~\ref{Lipschitztime} and~\eqref{parLipMsanseta}, it follows
\begin{eqnarray*}
2\epsilon_L &\leq & v (\bar{x}^*,t^*)- v (\bar{y}^*,t^*)\\
&\leq& v (\bar{x}^*,t^*)- v (\bar{x}^*,\hat{t}) + v (\bar{x}^*,\hat{t}) - v (\bar{y}^*,\hat{t})
+ v (\bar{y}^*,\hat{t}) -  v (\bar{y}^*,t^*)\\
&\leq& \omega_{\hat{t}}(|\bar{x}^*-\bar{y}^*|) + 2 \Lambda_0 |t^*-\hat{t}|,
\end{eqnarray*}
where $\omega_{\hat{t}}(\cdot)$ is the modulus of continuity associated with $v(\cdot,\hat{t}).$
Therefore, for all $\beta$ small enough, we get
\begin{eqnarray*}
\epsilon_L \leq  \omega_{\hat{t}}(|\bar{x}^*-\bar{y}^*|),
\end{eqnarray*}
which gives~\eqref{x-ybelow} with the modulus $\omega_{\hat{t}}$ independent
of $\beta$ yielding a lower bound for $|\bar{x}^* - \bar{y}^*|$, uniform in terms 
of $\beta$.

From this point, we continue the proof 
with the arguments given in the proof of Theorem~\ref{teo1}.
We point out that the constant $L$, which gives the
Lipschitz bound, depends on the constants appearing
in the assumptions, on the Lipschitz constant $L_0$ of $u_0$
and on $\mathrm{osc}_T (u)$ but not directly
neither on $||u(\cdot,t)||_\infty$ nor on $T.$

In particular, $L$ is independent on $||u_0||_{C^2(\R^d)}$
so we get the result for any Lipschitz continuous function $u_0$ by approximation.
\qed




\section{Application: large time behavior in the periodic setting.}
\label{LTBsection}

In this section we provide the large time behavior result for the problem~\eqref{pareqIto}-\eqref{initialdata} 
in the case when the datas are $\Z^d-$periodic. Hence, we argue on the problem
\begin{eqnarray}
\label{pareqtorus} 
&& \partial_t u - \mathrm{Tr}(A(x) D^2u) - \I^j(u(\cdot, t), x) + H(x, Du) = 0
\\
&& \nonumber \hspace*{7cm}  \mbox{in} \ \Q := \T^d \times (0,+\infty), \\
&& \label{initialtorus} u(\cdot,0)= u_0 \quad  \mbox{in} \ \T^d,
\end{eqnarray}
which is~\eqref{pareqIto}-\eqref{initialdata} in the periodic
setting with
\begin{equation*}\label{operator-periodic}
\I^j(u(\cdot, t), x) = \int_{\T^d} [u(x + j(x,z),t) - u(x,t) - \mathbf{1}_B (z)\langle Du(x,t), j(x, z)\rangle ] \nu(dz).
\end{equation*}

In order to have comparison principle/well-posedness results of Section~\ref{comparisonsection} and the 
regularity results given in 
Section~\ref{Lipschitzsection}, from now on, we assume that (A), (M'), (H0)-(H1)-(H2'), (J1)-(J2) hold with periodic datas
with respect to $x,z.$


\subsection{Solvability of the ergodic problem.}

\begin{prop}\label{propergodic}
There exists a unique constant $c \in \R$ for which the stationary ergodic problem
\begin{equation}\label{ergodic}
-\mathrm{Tr}(A(x)D^2 u) - \I^j(u, x) + H(x,Du) = - c, \quad \mbox{in} \ \T^d 
\end{equation}
has a solution $w \in W^{1,\infty}(\T^d)$. 
\end{prop}

In the proof of the above proposition we require an appropriate compactness property over the family of 
solutions $\{ u_\lambda \}$ of problem~\eqref{eq} (in the torus) as $\lambda \to 0$. This is the purpose
of the following lemma, whose proof follows closely the arguments of~\cite{Ley-Vihn}.


\begin{lema}\label{lemaosc}
Let $\lambda > 0$ and let $u$ be a continuous solution to~\eqref{eq} in $\T^d$. Then, there exists $C >0$
not depending on $\lambda$ such that $\mathrm{osc}(u) \leq C$.
\end{lema}


\noindent
{\bf \textit{Proof:}} We start claiming that, under Assumption (H1), 
there exists a constant $\eta > 0$ just depending on the data, 
and a sequence $L\to +\infty$ such that $H$ satisfies
\begin{equation}\label{Hcoercive}
H(x, L p) - L H(x, p) \geq \eta L^{m} |p|^m - \eta^{-1}
\end{equation}
for $x, p \in \R^d$. Note that, in particular, this proves that $H$ is superlinear.
Now, we consider $L \geq 1$ to be fixed and 
\begin{equation*}
M := \max \limits_{x, y \in \T^d} \{ u(x) - Lu(y) + (L - 1)\min \{ u \} - L |x - y|\}.
\end{equation*}

Note that if there exists $L$ such that $M \leq 0$, 
then for all $x, y \in \T^d$ we can write
\begin{equation*}
u(x) - Lu(y) \leq (1 - L)\min \{ u \} + L \sqrt{d},
\end{equation*}
and hence, taking $x, y \in \T^d$ such that $u(x) = \max \{ u \}$ and $u(y) = \min \{ u \}$, we get the result
with $C =L \sqrt{d}$.

Then, we assume that $M > 0$ for all $L\geq 1.$ 
In particular, the maximum in $M$ is attained at $(\bar x, \bar y)$ with $\bar x \neq \bar y$.
Thus, denoting $\varphi(x,y) = -(L-1)\min \{ u \} + L|x - y|$,
we can use $x \mapsto \varphi(x, \bar y)$ as a test function for $u$ at $\bar x$ and 
$y \mapsto -\varphi(\bar x, y)$ as a test function for $v := Lu$ at $\bar y$. Then, for each $\delta > 0$ we see that
\begin{equation*}
\begin{split}
\lambda u(\bar x) - \Tr(A(\bar x)X) - \I^j[B_\delta](\varphi(\cdot, \bar y), \bar x) 
- \I^j[B_\delta^c](u, \bar p, \bar x) + H(\bar x, L \bar{p}) \leq & \ 0 \\
\lambda v(\bar x) - \Tr(A(\bar y)Y) - \I^j[B_\delta](\varphi(\bar x, \cdot), \bar y) 
- \I^j[B_\delta^c](v, \bar p, \bar y) + L H(\bar x, \bar{p}) \geq & \ 0,
\end{split}
\end{equation*}
where $\bar p = (\bar x - \bar y)/ |\bar x - \bar y|$ and $X, Y$ satisfy~\eqref{matrixineq} with $\varphi = \varphi_1$. 

We subtract both inequalities and estimate each term arising in this operation. Note that by~\eqref{uboundstationary} we have
\begin{equation*}
\lambda (u(\bar x) - v(\bar y)) \leq (1 + L) H_0,
\end{equation*}
by Lemma~\ref{lemaADvarphi} we have
\begin{equation*}
\Tr(A(\bar x)X) - \Tr(A(\bar y)Y) \leq CL |\bar x - \bar y|,
\end{equation*}
for some constant $C > 0$ just depending on the data. Considering $\delta$ small in terms 
of $|\bar x - \bar y|$,
we use~\eqref{calcul-deriv134} (M') and (J1) similarly as in~\eqref{BdeltaLip} to get
\begin{equation*}
\I^j[B_\delta](\varphi(\cdot, \bar y), \bar x) 
 + \I^j[B_\delta](\varphi(\bar x, \cdot), \bar y) \leq L |\bar x - \bar y|^{-1}o_\delta(1).
\end{equation*}

By using that $(\bar x, \bar y)$ is the maximum point in $M$ we obtain
\begin{equation*}
\begin{split}
& \I^j[B_\delta^c](u, \bar p, \bar x) 
- \I^j[B_\delta^c](v, \bar p, \bar y) \\
\leq & \ L \int_{B_\delta^c} (|\bar x - \bar y + j(\bar x, z) - j(\bar y, z)| - |\bar x - \bar y| 
- \mathbf{1}_B(z) \langle \bar p, j(\bar x, z) - j(\bar y, z)\rangle) \nu(dz),
\end{split}
\end{equation*}
and performing a similar analysis as the one done in the proof of 
Theorem~\ref{teo1} (see~\eqref{proche0}-\eqref{plusloin0}),
we conclude that
\begin{equation*}
\I^j[B_\delta^c](u, \bar p, \bar x) 
- \I^j[B_\delta^c](v, \bar p, \bar y) \leq CL |\bar x - \bar y|,
\end{equation*}
for $C > 0$ and all $L$ large in terms of the data, and not depending on $\delta$.

Finally, using~\eqref{Hcoercive} we see that
\begin{equation*}
H(\bar x, L \bar p) - L H(\bar y, \bar p) \geq \eta L^m - \eta^{-1},
\end{equation*}
for $L$ large enough.

Then, joining the above estimates we conclude that there exists $C > 0$ such that, 
for all $L$ large just in terms of 
the data, we get
\begin{equation*}
\eta L^m \leq C(1+L) + L |\bar x - \bar y|^{-1} o_\delta(1),
\end{equation*}
and therefore, taking $\delta \to 0$ and enlarging $L$ if this is necessary, we arrive to a contradiction.
\qed

\medskip

\noindent
{\bf \textit{Proof of Proposition~\ref{propergodic}:}} 
For $\lambda > 0$, we denote $u_\lambda$ the unique bounded 
uniformly continuous solution to~\eqref{eq} 
given by Proposition~\ref{comparisonstationary}. 
By Lemma~\ref{lemaosc}, ${\rm osc}(u_\lambda)$ is bounded
independently of $\lambda.$ Therefore, by Theorem~\ref{teo1},
$u_\lambda$ is Lipschitz continuous with a constant independent of $\lambda.$

Defining 
\begin{equation*} 
v_\lambda(x) = u_\lambda(x) - u_\lambda(0), \quad x\in\T^d,
\end{equation*}
it follows that the family $\{ v_\lambda \}_{\lambda >0}$
is bounded and equi-Lipschitz continuous in $C(\T^d).$
By standard viscosity arguments, 
$v_\lambda$ satisfies the equation
\begin{equation*} 
\lambda v - \mathrm{Tr}(A(x) D^2v) - \I^j(v, x) + H(x, Dv) 
= -\lambda u_\lambda(0) \quad \mbox{in} \ \T^d,
\end{equation*}
where $\lambda u_\lambda(0)$ is bounded as $\lambda \to 0$ by~\eqref{uboundstationary}. 
Finally, by Ascoli Theorem and stability, up to subsequences,
there exists $w \in W^{1,\infty}(\T^d)$ and $c\in\R$
such that $v_\lambda \to w,$ $\lambda u_\lambda(0) \to c$ as $\lambda \to 0$ and the pair $(w,c)$ 
satisfies~\eqref{ergodic}.

For the uniqueness of $c$, we note that if there exist two pairs $(w_i, c_i), i=1,2$, both solutions to~\eqref{ergodic}, then the functions
$(x,t) \mapsto w_i(x) + c_i t, i=1,2$ are bounded viscosity solutions to equation~\eqref{pareqIto}. Hence, comparing them by the use of 
Proposition~\ref{parcomparison}, we obtain
\begin{equation*}
(c_1 - c_2)t \leq 2 ||w_1 - w_2||_\infty, \quad \mbox{for all} \ t > 0.
\end{equation*}

Dividing by $t$ and letting $t \to \infty$ we conclude $c_1 \leq c_2$. Since we can exchange the roles of $c_1$ and $c_2$, we conclude
the uniqueness of $c$.
\qed




\subsection{Strong maximum principle.} 

Before presenting our strong maximum principle result, we need to introduce some notation. 
For a L\'evy measure $\nu$ and a jump function $j$, 
for each $x \in \T^d$ we define  the push-forward
measure $\nu_x^j$ associated to $\nu$ through the 
function $z \mapsto j(x, z)$, that is, for each borel set $A \subset \R^d$
we have $\nu_x^j(A) = \nu(j^{-1}(x, A))$.
Thus, for each $x \in \T^d$ we define
\begin{equation*}
X_0(x) = \{x\}, \quad X_{n + 1} = \bigcup_{\xi \in X_n} \{ \xi + \mathrm{supp} \{ \nu_\xi^j \} \}, 
\ n \in \N,
\end{equation*}
where $\mathrm{supp}$ denotes the support of the measure, and the set
\begin{equation*}
\mathcal{X}(x) = \overline{\bigcup_{n \in \N} X_n(x)}. 
\end{equation*}

Finally, for $x \in \R^d$ we denote $E_0(x)$ the eigenspace associated to the null eigenvalue of $A(x)$.


\begin{prop}\label{SMaxP}{\bf (Strong maximum principle)}
Assume the hypotheses of Theorem~\ref{teo1} hold, with in addition, $H(x,p)$ locally Lipschitz in $p$. Assume the existence of a constant $r_0 > 0$ such that, for each $x \in \T^d$, 
\begin{equation}\label{covering}
B_{r_0}(x) \cap \{ x + E_0(x) \} \subset \mathcal{X}(x).  
\end{equation}

Let $u, v \in C(\mathcal{Q})$ be two solutions to~\eqref{pareqtorus}, associated to Lipschitz initial datum $u_0, v_0$, respectively.
Assume that $u-v$ achieves a maximum in $\Q=\T^d\times (0,+\infty)$
at $(x_0, t_0),$ that is,
\begin{equation*}
(u - v)(x_0, t_0) = \sup \limits_{\Q} \{ u - v\}.
\end{equation*}

Then, the function $u - v$ is constant in $\mathbb{T}^d \times [0, t_0]$. Moreover, we have
\begin{equation*}
(u - v)(x,t) = \sup \limits_{x \in \T^d} \{ u_0(x)-v_0(x)\}, \quad \mbox{for all} \ (x,t) \in \bar{\Q}.
\end{equation*}
\end{prop}


We remark that by the available Lipschitz regularity results, it is possible to reduce Proposition~\ref{SMaxP}
to a linear framework and this is the aim of the following lemma. We would like to stress on the role of the 
assumption~\eqref{covering} in its proof: in the directions of the second-order uniform ellipticity of the matrix $A$
the propagation of maxima follows classical arguments, and therefore the mixed operator extend this propagation
in the directions of degeneracy of $A$ through the sequential covering property of the nonlocal operator.
\begin{lema}\label{SMaxPlinear}
Let $A$ be a matrix satisfying (A), $\I^j$ as in~\eqref{operator} with $\nu$ satisfying (M') and
$j$ satisfying (J1), and both $\nu, j$ satisfying assumption~\eqref{covering}. 
Let $\beta\in L^\infty(\Q; \R^d)$ and $w$ be a bounded USC viscosity subsolution 
to the problem
\begin{equation}\label{pareqlinear}
\partial_t w - \mathrm{Tr}(A(x)D^2w) - \I^j(w, x) + \langle \beta(x,t), Dw\rangle = 0 \quad \mbox{in} \ \Q.
\end{equation}

If there exists $(x_0, t_0) \in \Q$ such that $M := w(x_0, t_0) = \sup \limits_{\Q} \{ w \}$, then $w(x,t_0) = M$ for all $x \in \T^d$.
\end{lema}


\noindent
{\bf \textit{Proof:}} It is sufficient to prove that under the assumptions of the problem, 
if $(x_0, t_0)$ is a global maximum point for $w$ subsolution to~\eqref{pareqlinear}, 
then $w$ is constant equal to $w(x_0, t_0)$ in 
$B_{r_0/4}(x_0) \times \{ t_0 \}$ where $r_0$ appears in~\eqref{covering}.

In fact, iterating the argument presented below a finite number of times, we conclude the main result.
Denote 
$$
K = \{ x \in \T^d : w(x, t_0) = M \}
$$ 
which is a nonempty closed set containing $x_0$.

Let $x^*\in K.$ Since $(x^*, t_0)$ is a global maximum point for $w$ we can use a constant function as a 
test function for $w$ at $(x^*, t_0)$. Thus, for each $\delta > 0$ we have
\begin{equation*}
-\int_{B_\delta^c} [w(x^* + j(x^*, z), t_0) - w(x^*, t_0)]\nu(dz) \leq 0,
\end{equation*}
and since $w(x^* + j(x^*, z), t_0) \leq w(x^*, t_0),$ we obtain
\begin{equation*}
\int_{B_\delta^c} [w(x^* + j(x^*, z), t_0) - w(x^*, t_0)]\nu(dz) = 0,
\end{equation*}

Therefore, since $z\mapsto w(x^*+j(x^*,z),t_0)-w(x^*,t_0)$
is upper semicontinuous and $\delta > 0$ is arbitrary, we get $w(x,t_0) = M$
for each $x \in x^* + \mathrm{supp} \{ \nu_{x^*}^j \}$. We apply the same argument inductively 
to conclude that $\mathcal{X}(x^*) \subseteq K$.

Noting that $\mathcal{X}(x) \subseteq \mathcal{X}(x^*)$ for all $x \in \mathcal{X}(x^*)$, 
by the use of~\eqref{covering} we have
\begin{equation}\label{keycovering}
B_{r_0}(x) \cap \{ x + E_0(x) \} \subset \mathcal{X}(x^*), \quad 
\mbox{for each $x \in \mathcal{X}(x^*)$ and $x^*\in K.$} 
\end{equation}

Consider the open set 
$\Gamma := B_{r_0/4}(x_0) \setminus K$. If $\Gamma = \emptyset$, then the result follows. 
From now on, we argue by contradiction assuming that $\Gamma \not= \emptyset.$ It follows that
there exists $\bar{x} \in \Gamma$,
$0 < R < r_0/4$ and $x^* \in \partial K$  such that
\begin{equation*}
B_R(\bar{x}) \subset \Gamma \quad \mbox{and} 
\quad  x^*\in \partial B_R(\bar{x}) \cap K.
\end{equation*}

Up to replace $\bar{x}$ by $(\bar{x}+x^*)/2$ if needed, we may assume $\partial B_R(\bar{x}) \cap K=\{x^*\}.$

At this point, for $\gamma, h > 0$ to be fixed, we introduce the function 
\begin{equation*}
\phi(x,t) = e^{-\gamma R^2} - e^{-\gamma d(x,t)}
\quad \mbox{with $d(x,t) = |x - \bar{x}|^2 + h(t - t_0)^2.$}
\end{equation*}

Direct computations say that for each $(x,t) \in \T^d$
\begin{equation*}
\begin{split}
\partial_t \phi(x,t) & = 2\gamma h e^{-\gamma d(x,t)} (t - t_0) \\
D\phi(x,t) & = 2\gamma e^{-\gamma d(x,t)} (x - \bar{x}) \\
D^2\phi(x,t) & = 2\gamma e^{-\gamma d(x,t)} [I_d - 2\gamma (x - \bar{x}) \otimes (x - \bar{x})],
\end{split}
\end{equation*}
meanwhile, following~\cite{Ciomaga}, there exists $C_{\nu, j} > 0$ depending only on the data, such that
\begin{equation*}
\I^j(\phi(\cdot, t), x) \leq \gamma e^{-\gamma d(x,t)} C_{\nu, j} .
\end{equation*}

With these estimates and applying (A), for each $x \in \T^d$ we have 
\begin{eqnarray*}
&& \mathcal{E}(\phi, x, t_0) \\
& := &  \partial_t \phi(x, t_0) - \mathrm{Tr}(A(x)D^2\phi (x, t_0)) 
- \I^j(\phi(\cdot, t_0), x) + \langle \beta(x,t), D\phi(x, t_0)\rangle \\
& \geq & 2\gamma e^{-\gamma d(x,t)} \Big{(} h(t - t_0) - L_\sigma + \gamma |\sigma^T(x) (x - \bar{x})|^2
- C_{\nu, j} -  ||\beta||_\infty |x - \bar{x}| \Big{)}
\end{eqnarray*}

Note that $\bar{x} \in B_{r_0/4}(x^*)$. If $\bar{x} - x^* \in E_0(x^*)$, by~\eqref{keycovering} we would have 
$\bar{x} \in \mathcal{X}(x^*) \subseteq K$, which is a contradiction with the choice of $\bar{x}$. 
Thus, $\sigma^T(x^*) (x^* - \bar{x})\not=0.$ By continuity, there exists $\eta(x^*), R^* > 0$
such that
\begin{equation*}
|\sigma^T(x) (x - \bar{x})| \geq \eta(x^*) \quad \mbox{for all} \ x \in B_{R^*}(x^*).
\end{equation*}

This allows us to get the inequality
\begin{equation*}
\mathcal{E}(\phi, x, t_0) 
\geq  2\gamma e^{-\gamma d(x,t)} \Big{(} h|t - t_0| + \gamma \eta(x^*)^2 -L_\sigma
- C_{\nu, j} -  2R ||\beta||_\infty  \Big{)},
\end{equation*}
for each $x \in B_{R^*}(x^*)$. Thus, taking $\gamma$ large in terms of $R, h, t_0$ and the data, 
we conclude $v$ that is a strict supersolution to~\eqref{pareqlinear} 
in $B_{R^*}(x^*) \times (0, t_0 + 1)$. 

On the other hand, since $\bar{B}_R(\bar{x}) \cap K=\{x^*\},$
there exists $\rho^* >0$ such that
\begin{equation*}
w(x, t_0) \leq M - \rho^* \quad \mbox{for all} \ x \in \bar{B}_{R}(\bar{x}) \setminus B_{R^*}(x^*),
\end{equation*}
and therefore, by upper semicontinuity of $w$, there exists $\tau^* \in (0,1)$ small enough such that
\begin{equation}\label{w<Mint}
w \leq M - \rho^*/2 \quad \mbox{in} \
(\bar{B}_{R}(\bar{x}) \setminus B_{R^*}(x^*)) \times (t_0 -\tau^*, t_0 + \tau^*).
\end{equation}

At this point, we fix $h > (R / \tau^*)^2$. Under this choice, the ellipsoid 
$$
\varSigma = \{ (x,t) : |x - \bar{x}|^2 + h(t - t_0)^2 \leq R^2 \}
$$
satisfies $\varSigma \subset \bar{B}_R(\bar{x}) \times (t_0 - \tau^*, t_0 + \tau^*)$. 
Notice that $(x^*,t_0)\in\partial\varSigma$ since $|x^*-\bar{x}|=R$
and  $w(x^*,t_0)=M.$
Since $\phi > 0$ in $\varSigma^c$, 
for all $\epsilon > 0$ we have $w - \epsilon \phi < M$ in $\varSigma^c$, and by~\eqref{w<Mint}, 
taking $\epsilon > 0$ small in terms of $\rho^*$, we obtain
\begin{eqnarray*}
w-\epsilon\phi \leq M - \rho^*/2 +\epsilon ||\phi||_{L^\infty(\varSigma)}< M
\ \mbox{in $\varSigma \setminus (B_{R^*}(x^*) \times (t_0 - \tau^*,t_0 + \tau^*)).$}
\end{eqnarray*}

Hence, we conclude from this that $w - \epsilon \phi$ attains its global maximum at a point 
$(x', t') \in \varSigma$ with $x' \in B_{R^*}(x^*)$.
Since $w$ is a viscosity subsolution to~\eqref{pareqlinear}, we get
$
\mathcal{E}(\epsilon \phi, x', t') \leq 0.
$

By the linearity of~\eqref{pareqlinear} this drives us to the inequality
\begin{equation*}
\mathcal{E}(\phi, x', t') \leq 0, 
\end{equation*}
which contradicts the fact that $v$ is a strict supersolution to~\eqref{pareqlinear}
in $B_{R^*}(x^*) \times (0, t_0 + 1).$
\qed

\medskip

The following lemma is a consequence of the comparison principle, see~\cite{Barles-Souganidis}.
\begin{lema}\label{lemaSMP}
Assume assumptions of Proposition~\ref{parcomparison} hold. Let $u, v$ be respectively a
bounded USC subsolution and a bounded LSC supersolution to equation~\eqref{pareqtorus}
and for $t \in [0,+\infty)$, define 
\begin{equation*} 
\kappa(t) = \sup \limits_{x \in \T^d} \{ u(x, t) - v(x, t)\}. 
\end{equation*}

Then, for all $0 \leq s \leq t$, we have $\kappa(t) \leq \kappa(s)$.
\end{lema}


The previous lemmas allows to provide the

\medskip

\noindent
{\bf \textit{Proof of Propostion~\ref{SMaxP}:}} By Lemma~\ref{lemaSMP}, 
the continuity of $u-v$ and the fact that $(x_0, t_0)$ 
is a global maximum point for $u - v$, we have 
$(u-v)(x_0,t_0)=\kappa(0) = \kappa(\tau)$ for all $\tau \in [0,t_0]$. 
Then, it is sufficient to prove that
for each $\tau \in (0, t_0]$, $(u - v)(x, \tau) = \kappa(0)$ for all $x \in \T^d$, concluding the result up to $\tau = 0$ by continuity.

By Proposition~\ref{Lipspaceparabolic}, $u$ and $v$ are Lipschitz in space in $[0,t_0]$, 
with Lipschitz constant depending only 
on the data and $t_0.$ Then, by classical arguments in the viscosity theory, 
the function $w := u - v$ is a viscosity subsolution to the problem
\begin{equation*}
\partial w - \mathrm{Tr}(A(x)D^2w) - \I^j(w(\cdot, t), x) + \langle \beta(x,t), Dw \rangle \leq 0 \quad \mbox{in} \ \T^d \times (0, t_0],
\end{equation*}
with $\beta \in L^\infty(\T^d \times [0, t_0]; \R^d)$ defined as
\begin{equation*}
\beta(x,t) = \int_{0}^{1} D_p H(x, sDu(x, t) + (1 - s)Dv(x, t)) ds.
\end{equation*}

Therefore, for all $\tau\in (0,t_0],$ there exists
$x_\tau\in\T^d$ such that $w(x_\tau,\tau)=\kappa(\tau).$
By Lemma~\ref{SMaxPlinear}, we obtain
$w(\cdot,\tau)=\kappa(\tau)=\kappa(0)$ and the results follows.
\qed


\subsection{Large time behavior.} 
The above results are sufficient to get the large time behavior 
for~\eqref{pareqtorus}-\eqref{initialtorus}. 
\begin{prop}\label{propuniquenessergodic}
Under the assumptions of this section and of Proposition~\ref{SMaxP}, 
the continuous solution of~\eqref{ergodic} is unique up to a constant.
\end{prop}

\begin{teo}\label{teoLTB}{\bf (Ergodic large time behavior)}
Under the assumptions of this section, for any $u_0\in W^{1,\infty},$ there
exists a unique solution $u$ to problem~\eqref{pareqtorus}-\eqref{initialtorus}. 
Under the additional assumptions of Proposition~\ref{SMaxP}, 
there exists a pair $(u_\infty, c)$ solution 
to~\eqref{ergodic} such that, as $t \to \infty$
\begin{equation*}
u(\cdot, t) + c t \to u_\infty \quad \mbox{in} \ W^{1,\infty}(\T^d).
\end{equation*}
\end{teo}

The proof of Proposition~\ref{propuniquenessergodic} is an easy consequence
of Propositions~\ref{parcomparison} and~\ref{SMaxP} and follows the same lines as 
in~\cite{Barles-Chasseigne-Ciomaga-Imbert, Barles-Koike-Ley-Topp}.
To prove Theorem~\ref{teoLTB}, we first notice that, by comparison,
for every solution $(v,c)\in W^{1,\infty}(\T^d)\times \R$ 
of~\eqref{ergodic}, there exists $M>0$
such that $v(x)-M\leq u(x,t)+ct \leq v(x)+M.$ It follows that
${\rm osc}(u(\cdot,t))$ is bounded independently of $t.$ 

Therefore, by Proposition~\ref{Lipschitztime}, $u(\cdot,t)$ is Lipschitz
continuous with a constant independent of $t.$ Hence,
$\{u(\cdot,t)+ct, t\geq 0\}$ is relatively compact
in $W^{1,\infty}(\T^d).$ The proof of the convergence of the whole
sequence then follows as in~\cite{Barles-Chasseigne-Ciomaga-Imbert}.

\appendix
\section{}
\label{secappendix}

We provide the technical estimates used in Theorem~\ref{teo1}. 

We start with the following relationship between $\I^j$ and $\J^j$, see~\cite{Ciomaga} for a proof. For $\J^j$ defined in~\eqref{defJ} we adopt the analogous notations as those introduced for $\I^j$
at the end of the introduction.
\begin{lema}\label{lemaIexp}
Let $g \in C^2(\R^d)$. Then, for each $\delta \in (0,1)$ we have
\begin{equation*}
\J^j[B_\delta](g, x) = \I^j[B_\delta](g, x) 
+ ||Dg||_{L^\infty(B_\delta(x))}^2 O(\delta^{2 - \sigma}),
\end{equation*}
where the $O$-term is independent of $g$.
\end{lema}

Next result is useful in the linearization of the exponential terms arising in $\J^j$.
\begin{lema}\label{enlt1}
Let $g: \R^d \times \R^d \to \R$ a bounded measurable function. Assume there exist $L > 1$, $b \in (0,1)$, $C_1, C_2 > 0$, $\mathcal P \subset \R^d$ measurable and $\bar x, \bar y \in \R^d$ with $|\bar x - \bar y| \leq 1$, such that for all $z \in \mathcal P$ we have
\begin{eqnarray}
\label{cota1} |g(\bar x, z) - g(\bar y, z)| \leq C_1(L |\bar x- \bar y| + b)|z|, \ \mbox{and} \\
\label{cota2} |g(\bar x, z)|, |g(\bar y, z)| \leq \min \{ C_2, C_1 (L + b) |z|\}.
\end{eqnarray}

Then, there exists $C_3$ just depending on $C_1, C_2$ such that
\begin{equation*}
e^{g(\bar x, z)} - e^{g(\bar y, z)} \leq g(\bar x, z) - g(\bar y, z)
+ C_3 (L^2|\bar x -\bar y| + Lb)|z|^2, \quad \mbox{for} \ z \in \mathcal P.
\end{equation*}
\end{lema}


\noindent
{\bf \textit{Proof:}}
By the Mean Value Theorem, for each $z \in \mathcal P$ we have
\begin{equation*}
\begin{split}
e^{g(\bar x, z)} - e^{g(\bar y, z)} = & \ e^{g(\bar y, z)} (e^{g(\bar x, z) - g(\bar y, z)} - 1) \\
\leq & \ e^{g(\bar y, z)} (g(\bar x, z) - g(\bar y, z) + e^{\xi_1(z)}(g(\bar x, z) - g(\bar y, z))^2/2) 
\end{split}
\end{equation*}
for some $\xi_1(z) \in [-2C_2, 2 C_2]$ in view of~\eqref{cota2}. Using this and again the Mean Value Theorem (this time on the term $e^{g(\bar y, z)}$) we obtain
\begin{equation*}
\begin{split}
e^{g(\bar x, z)} - e^{g(\bar y, z)} \leq & \ g(\bar x, z) - g(\bar y, z) 
+ (e^{g(\bar y, z)} - 1)(g(\bar x, z) - g(\bar y, z)) \\
& + \frac{e^{3C_2}}{2} |g(\bar x, z) - g(\bar y, z)|^2 \\
\leq & \ g(\bar x, z) - g(\bar y, z) + e^{\xi_2(z)}|g(\bar y, z)||g(\bar x, z) - g(\bar y, z) | \\
& + \frac{e^{3C_2}}{2} |g(\bar x, z) - g(\bar y, z)|^2
\end{split}
\end{equation*}
with $\xi_2(z) \in [-C_2, C_2]$. Then, using this last fact together with~\eqref{cota1} and~\eqref{cota2} we infer
\begin{equation*}
\begin{split}
e^{g(\bar x, z)} - e^{g(\bar y, z)} \leq & \ g(\bar x, z) - g(\bar y, z)
+ C_1^2 e^{C_2} (L + b)(L|\bar x -\bar y| + b)|z|^2 \\
& + \frac{e^{3C_3}}{2} C_1^2(L |\bar x - \bar y| + b)^2|z|^2 \\
\leq & \ g(\bar x, z) - g(\bar y, z) + C |z|^2 (L|\bar x -\bar y| + b) (L
+ L |\bar x - \bar y| + b)
\end{split}
\end{equation*}
for some $C > 0$ just depending on $C_1$ and $C_2$. Since $b, |\bar x -\bar y| \leq 1 \leq L$ we conclude the proof.
\qed


\begin{lema}\label{estim-nonlocale1}
Let $\bar x, \bar y \in \R^d$ with $|\bar x - \bar y| > 0$ and denote $\bar p = (\bar x-\bar y)/|\bar x- \bar y|$. Define, for $z \in \R^d$ the function 
\begin{equation*}
\Psi(z) = |\bar x- \bar y + j(\bar x , z) - j(\bar y, z)| - |\bar x - \bar y| - \langle \bar p, j(\bar x, z) - j(\bar y, z)\rangle.
\end{equation*}

Then, there exists $C > 0$ just depending on the data such that, for all $\delta > 0$ and $\mathcal P \subset \R^d$ measurable, we have the estimate
\begin{equation*}
\int_{\mathcal P \cap B \setminus B_\delta} \Psi(z) \nu(z)dz \leq C |\bar x - \bar y|.
\end{equation*}
\end{lema}

\noindent
{\bf \textit{Proof:}} Getting nonnegative upper bounds for $\Psi$ in the domain of integration, we can get rid of the intersection with $\mathcal P$ and therefore we omit it for simplicity. 

Notice that for $z$ such that $C_j |z| \leq 1/2$, by (J1) we have
\begin{equation*}
|j(\bar{x},z) - j(\bar{y},z)| \leq |\bar{x} - \bar{y}|/2.
\end{equation*}

Then, for $|z| \leq (2C_j)^{-1}$ we can perform a Taylor expansion on $\Psi(z)$ (around the point $\bar x- \bar y $) to write
\begin{equation} \label{proche0}
\Psi(z) \leq \frac{2}{|\bar{x} - \bar{y}|}|j(\bar{x},z) - j(\bar{y},z)|^2
\leq  2 C_j^2 |\bar{x} - \bar{y}||z|^2\; 
\end{equation}
and introducing $\delta_0 = (2C_j)^{-1}$, by the above inequality and (M') we get
\begin{equation*}
\int_{B_{\delta_0} \setminus  B_\delta} \Psi(z) \nu(dz) \leq C L |\bar x - \bar y|,
\end{equation*}
where $C > 0$ just depend on the data. 

On the other hand, when $|z| > (2C_j)^{-1}$ we apply triangular inequality, and the fact that $|\bar p|=1$ 
together with Cauchy-Schwartz inequality to get, by using (J1), the inequality
\begin{equation} \label{plusloin0}
\Psi(z) \leq  2|j(\bar{x}, z) - j(\bar{y}, z)|
\leq  2 C_j |\bar{x} - \bar{y}||z|.
\end{equation}

Using this and (M') we can write
\begin{equation*}
\int_{B \setminus B_{\delta_0}} \Psi(z) \nu(dz) \leq C|\bar x - \bar y|,
\end{equation*}
for some $C > 0$ depending on the data. This concludes the proof.
\qed

\bigskip

\noindent {\bf Acknowledgements.} 
E. T. was partially supported by Fondecyt Postdoctoral Grant No. 3150100 and Conicyt PIA Grant No. 79150056. This work was partially supported by the ANR (Agence Nationale de
la Recherche) through projects HJnet ANR-12-BS01-0008-01 and WKBHJ ANR-12-BS01-0020.



\begin{thebibliography}{00}
\bibitem{Alvarez-Tourin}
Alvarez, O and Tourin, A. {\em Viscosity Solutions of Nonlinear Integro-Differential Equations} 
Ann. Inst. H. Poincar\'e Anal. Non Lin\'eaire 13 (3) (1996) 293-317.

\bibitem{Bardi-DaLio}
Bardi, M. and Da Lio, F. {\em On the Strong Maximum Principle for Fully Nonlinear Degenerate Elliptic Equations} Arch. Math (Basel), 73 (4),
276-285, 1999.



\bibitem{Barles-book}
Barles, G. {\em Solutions de Viscosit\'e des \'Equations de Hamilton-Jacobi} Collection ``Math\'ematiques et Applications'' de la SMAI, no 17, Springer-Verlag (1994).

%

\bibitem{Barles1}
Barles, G. {\em A Short Proof if the $C^{0,\sigma}-$regularity of Viscosity Subsolutions for Superquadratic Viscous Hamilton-Jacobi Equations and Applications} Nonlinear Analysis 73 (2010) 31-47.



\bibitem{Barles2}
Barles, G. {\em A Weak Bernstein Method for Fully Nonlinear Elliptic Equations.} Diff. and Integral Equations, 4(2): 241-262, 1991.

\bibitem{BCI}
Barles, G., Chasseigne, E. and Imbert, C., {\em H\"older continuity of
  solutions of second-order non-linear elliptic integro-differential
  equations}, J. Eur. Math. Soc. (JEMS), 13 (2011), pp.~1--26.


\bibitem{Barles-Chasseigne-Imbert-Ciomaga-lip}
Barles, G., Chasseigne, E., Ciomaga, A. and Imbert, C. {\em Lipschitz Regularity of Solutions for Mixed Integro-Differential Equations.}
J. Diff. Eq., 252 (2012), 6012-6060.



\bibitem{Barles-Chasseigne-Ciomaga-Imbert}
Barles, G., Chasseigne, E., Ciomaga, A. and Imbert, C. {\em Large time behavior of periodic viscosity solutions for uniformly parabolic integro-differential equations.} 
Calc. Var. Partial Differential Equations 50 (2014), 283--304.






%



%

\bibitem{Barles-Imbert}
Barles, G. and Imbert, C. {\em Second-order Eliptic Integro-Differential Equations: Viscosity Solutions' Theory Revisited.} Ann. Inst. H. Poincar\'e Anal. Non Lin\'eaire 25 (2008) 567-585.



\bibitem{Barles-Koike-Ley-Topp}
Barles, G., Koike, S., Ley,  O. and Topp, E. {\em Regularity results and large time behavior for integro-differential equations with coercive Hamiltonians.} Calc. Var. Partial Differential Equations 54 (2015), 539--572.

\bibitem{BP1}
Barles, G.  and Perthame, B. {\em Exit time problems in optimal
control and vanishing viscosity method.} SIAM J. in Control and
Optimization, {\bf 26}, 1988, pp~1133-1148.




\bibitem{Barles-Souganidis1}
Barles, G. and Souganidis, P.E. {\em On the large time behavior of solutions of Hamilton-Jacobi equations} SIAM J. Math. Anal. 31(4), 925–939 (2000)



\bibitem{Barles-Souganidis}
Barles, G. and Souganidis, P.E. {\em Space-time Periodic Solutions and Long-Time Behavior of Solutions of Quasilinear Parabolic Equations} SIAM J. Math. Anal., 32 (2001), 1311-1323 (electronic).



\bibitem{Barles-Topp}
Barles, G. and Topp, E. {\em Existence, Uniqueness and Asymptotic Behavior for Nonlocal Parabolic Problems with Dominating Gradient Terms.} Preprint.



\bibitem{Barles-Topp-censored}
Barles, G. and Topp, E. {\em Lipschitz Regularity for Censored Subdiffusive Integro-Differential 
Equations with Superfractional Gradient Terms.} Nonlinear Analysis 131 (2016), 3-31.



\bibitem{Bernstein}
Bernstein, S. {\em Sur la g\'en\'eralisation du probl\'eme de Dirichlet, I Math. Ann., 62 (1906) 253-271, II, Math. Ann., 69 (1910), 82-136}.






%

%




%

\bibitem{Capuzzo-Dolcetta-Leoni-Porretta}
Capuzzo-Dolcetta, I., Leoni, F. and Porretta, A. {\em H\"older Estimates for Degenerate Elliptic Equations with Coercive Hamiltonians.} Trans. Amer. Math. Soc. 362 (9) 4511-4536 (2010).



\bibitem{Ciomaga}
Ciomaga, A. {\em On the Strong Maximum Principle for Second Order Nonlinear Parabolic Integro-Differential Equations} Advances in Diff. Equations. 17 (2012), 635-671.



\bibitem{Coville}
Coville, J. {\em Maximum principles, sliding techniques and applications to nonlocal equations}
Electronic Journal of Differential Equations, Vol. 2007(2007), No. 68, pp. 1-23.



\bibitem{Coville2}
Coville, J. {\em Remarks on the Strong Maximum Principle for Nonlocal Operators} 
Electron. J. Differential Equations, No. 66 (2008) pp. 1-10.



\bibitem{usersguide} 
Crandall, M.G., Ishii H. and Lions, P.-L. {\em User's Guide to Viscosity Solutions of Second Order Partial Differential Equations.} Bull. Amer. Math. Soc. (N.S.), Vol. 27 (1992), no. 1, 1-67.



\bibitem{DaLio}
Da Lio, F. {\em Comparison Results for Quasilinear Equations in Annular Domains and Applications.} Comm. Partial Diff. Equations, 27 (1 \& 2) 283-323 (2002).



\bibitem{Hitch}
Di Neza, E., Palatucci, G. and Valdinoci, E. {\em Hitchhiker's Guide to the Fractional Sobolev Spaces.}
Bull. Sci. Math., 136, (2012), no. 5, 521--573.



\bibitem{GKLV}
Galise, G., Koike, S., Ley, O. and Vitolo, A. {\em Entire Solutions of Fully Nonlinear Elliptic Equations with a Superlinear Gradient Term.} Preprint.









%



%

%

\bibitem{Fr}Frehse, J.
{\em On the regularity of solutions to elliptic differential inequalities}. Mathematical techniques of optimization, control and decision, pp. 91–109, Birkhäuser, Boston, Mass., 1981. 

\bibitem{G-T} D. Gilbarg and N.S. Trudinger, \emph{Elliptic partial differential equations of second order},  Springer-Verlag, Berlin 2001.

%



%

\bibitem{Ishii2}
Ishii, H. {\em Perron's Method for Hamilton-Jacobi Equations.} Duke Math. J. 55 (1987), 369-384.



\bibitem{Ishii3}
Ishii, H. {\em On Uniqueness and Existence of Viscosity Solutions of Fully Nonlinear Second-Order Elliptic PDEs} Comm. Pure Appl. Math., 42(1):15–45, 1989.



\bibitem{Ishii-Lions}
Ishii, H. and Lions, P.L. {\em Viscosity Solutions of Fully Nonlinear Second-Order Elliptic Partial Differential Equations } J. Differential Equations, 83(1) 26-78, 1990.



\bibitem{Jensen}
Jensen, R. {\em The Maximum Principle for Viscosity Solutions of Fully Nonlinear Second Order Partial Differential Equations} Arch. Rational Mech. Anal., 101(1):1–27, 1988.



\bibitem{Lasry-Lions}
Lasry, J.M. and Lions, P.L. {\em Nonlinear elliptic Equations with Singular Boundary Conditions and Stochastic Control with State Constraints.}
Math. Ann. 283, 583-630 (1989).



\bibitem{Ley-Vihn}
Ley, O. and Nguyen, V.D. {\em Gradient bounds for nonlinear degenerate parabolic equations
and application to large time behavior of systems.} Nonlinear Anal. 130 (2016), 76--101.

\bibitem{book} Lions P. L. {\em Generalized solutions of Hamilton-Jacobi
equations}, Pitman, Boston, 1982.





%




%

\bibitem{Tchamba}
Tchamba, T.T. {\em Large Time Behavior of Solutions of Viscous Hamilton-Jacobi Equations with Superquadratic Hamiltonian.} Asymptot. Anal. 66 (2010) 161-186.
\end{thebibliography}
\end{document}